\documentclass[a4paper,11pt]{amsart}\usepackage{amsfonts,amssymb,bbm}
\usepackage{enumerate,mathrsfs,xfrac}

\usepackage[utf8]{inputenc}

\def\C{\Bbb{C}}
\def\k{\Bbbk}
\def\K{\mathbbm{K}}
\def\bk{{\bar{\k}}}

\def\N{\Bbb{N}}\def\Q{\Bbb{Q}}\def\R{\Bbb{R}}

\def\di{\partial}

\def\bl{\langle}\def\br{\rangle}

\def\RpX{R^{\oplus p}_X}\def\RpXu{R^{\oplus p}_{X,u}}\def\RpXt{R^{\oplus p}_{X,t}}\def\RpY{R^{\oplus p}_Y}

\newcommand{\quot}[2]{{\footnotesize\left.\raisebox{1.2ex}{$#1$}\!\! \ensuremath\diagup \!\!\raisebox{-1.2ex}{$#2$}\right.}}
\newcommand{\quots}[2]{{\footnotesize\left.\raisebox{0.4ex}{$#1$}\! / \!\raisebox{-0.4ex}{$#2$}\right.}}

\def\tc{\tilde{c}}\def\tf{{\tilde{f}}}\def\tF{\tilde{F}}
 
\def\tr{{\tilde{r}}}\def\tR{{\tilde{R}}}
\def\tt{{\tilde{t}}}

\def\tu{\tilde{u}}\def\tx{{\tilde{x}}}\def\tX{{\tilde{X}}}

\def\ha{\hat{a}}\def\he{\hat{e}}\def\hg{{\hat{g}}}
\def\hR{{\widehat{R}}}

\def\ga{\gamma}\def\de{\delta}
\def\ep{\epsilon}
\def\La{\Lambda}

\def\cA{\mathscr A}\def\ca{\mathfrak a}
\def\cb{\mathfrak b}\def\cC{\mathscr C}
\def\cF{\mathcal F}\def\cf{{\frak f}}
\def\cG{\mathscr G}\def\chG{\widehat{\cG}}
\def\cK{{\mathscr K}\!}\def\cL{\mathscr L}\def\cR{\mathscr{R}}

\def\cm{{\frak m}}
\def\cU{\mathcal U}\def\cu{{\frak u}}

\def\um{{\underline{m}}}\def\un{{\underline{n}}}

\def\one{{1\hspace{-0.1cm}\rm I}}\def\zero{\mathbb{O}}

\newcommand{\ber}{\begin{array}{l}}\newcommand{\eer}{\end{array}}
\newcommand{\bpm}{\begin{pmatrix}}\newcommand{\epm}{\end{pmatrix}}
\newcommand{\bM}{\begin{matrix}}\newcommand{\eM}{\end{matrix}}
\newcommand{\bee}{\begin{enumerate}}\newcommand{\eee}{\end{enumerate}}
\newcommand{\bei}{\begin{itemize}}\newcommand{\eei}{\end{itemize}}

\def\wrt{with respect to }
\def\sset{\subset}\def\sseteq{\subseteq}\def\smin{\setminus}

\def\MapX{Maps\big(X,(\k^p,o)\big)}

\newtheorem{Lemma}{Lemma}[section]\newcommand{\bel}{\begin{Lemma}}\newcommand{\eel}{\end{Lemma}}
\newtheorem{Theorem}[Lemma]{Theorem}\newcommand{\bthe}{\begin{Theorem}}\newcommand{\ethe}{\end{Theorem}}
\newtheorem{Proposition}[Lemma]{Proposition}\newcommand{\bprop}{\begin{Proposition}}\newcommand{\eprop}{\end{Proposition}}
\newtheorem{Corollary}[Lemma]{Corollary}\newcommand{\bcor}{\begin{Corollary}}\newcommand{\ecor}{\end{Corollary}}
\newtheorem{Definition}[Lemma]{Definition}\newcommand{\bed}{\begin{Definition}}\newcommand{\eed}{\end{Definition}}
\newtheorem{Definition-Proposition}[Lemma]{Definition-Proposition}

\def\bpr{{\em\noindent Proof.\ }}
\newcommand{\epr}{{\hfill\ensuremath\blacksquare}}
\newtheorem{Remark}[Lemma]{Remark}\newcommand{\beR}{\begin{Remark}\rm}\newcommand{\eeR}{\end{Remark}}
\newtheorem{Example}[Lemma]{Example}\newcommand{\bex}{\begin{Example}\rm}\newcommand{\eex}{\end{Example}}
\newtheorem{Problem}[Lemma]{Problem}\newcommand{\bprob}{\begin{Problem}\rm}\newcommand{\eprob}{\end{Problem}}

\newcommand{\bet}{\begin{tabular}{cccccccc}}\newcommand{\eet}{\end{tabular}}
\newcommand{\beq}{\begin{equation}}\newcommand{\eeq}{\end{equation}}

\newcommand{\isom}[1]{\xrightarrow[\,\smash{\raisebox{1.15ex}{\ensuremath{\scriptstyle\sim}}}\,]{#1}}
\title[]{U\MakeLowercase{nfoldings of maps, the first results on stable maps, and results of   }
 \\M\MakeLowercase{ather-}Y\MakeLowercase{au/}G\MakeLowercase{affney-}H\MakeLowercase{auser type in arbitrary characteristic.}
 }
\author[]{D\MakeLowercase{mitry} K\MakeLowercase{erner}}
\address{Department of Mathematics, Ben Gurion University of the Negev, P.O.B. 653, Be'er Sheva 84105, Israel.}
\email{dmitry.kerner@gmail.com}
\date{\today\ \  filename: \jobname.tex}
\thanks{I was supported by the Israel Science Foundation (grant No.  1910/18)}

\begin{document}\maketitle
\begin{abstract}
Consider the (formal/analytic/algebraic) map-germs $\MapX$. Let $\cG$ be the group of right/contact/left-right transformations.
 I extend the following (classical) results from the real/complex-analytic case to the case of arbitrary field    $\k$.
 \bei
 \item
 A separable unfolding is locally trivial iff it is infinitesimally trivial.
 \item
 An unfolding is locally versal iff it is infinitesimally versal.
\item The criterion of factorization of map-germs in zero characteristic. When is the map $X\to (\k^p,o)$ $\cG$-equivalent
 to the composition $X\to \tX\to (\k^p,o)$ with $X\cong \tX\times(\k^r,o)$?
\item Criteria of  trivialization of unfoldings over affine base.
\item Fibration of $\cK$-orbits into $\cA$-orbits.
\item A map is locally stable iff it is infinitesimally stable.
\item Stable maps are unfodings of their genotypes.
\item Stable maps are determined by their local algebras.
\item Results of Mather-Yau/Scherk/Gaffney-Hauser type. How does the module $T^1_\cG f$, or related algebras, determine the $\cG$-equivalence type of $f$?
 \eei

\end{abstract}
\setcounter{secnumdepth}{6} \setcounter{tocdepth}{2}

\section{Introduction}\label{Sec.Introduction}
\subsection{}
 Let $\k\in \R,\C$ and consider $\k$-analytic map-germs $(\k^n,o)\to(\k^p,o)$. These are studied up to the right ($\cR$), contact ($\cK$),
  and left-right ($\cA$) equivalences. Among the cornerstones of Singularity Theory are the finite determinacy, the theory of unfoldings and the
   theory of stable maps. (See e.g. \cite{AGV}, \cite{AGLV}, \cite{Gr.Lo.Sh}, \cite{Martinet.1982}, \cite{Mond-Nuno}.) An additional line of research was around
    the question ``How is a map determined by its critical/singular/instability locus?"
     One way to answer this goes via the classical results of \cite{Mather-Yau}, \cite{Gaffney-Hauser}, \cite{Scherk}.

The classical approach relied heavily on vector fields integration. Initially numerous results could not be extended to the case ``$\k$ is any field", not even to
 the zero characteristic case. The study of $\cR,\cK$ equivalences over an arbitrary field began in \cite{Greuel-Kroning}. Many results on $\cR,\cK$
  are available
  by now, see e.g. \cite{BGK.20}, \cite{B.K.motor}, \cite{Greuel-Pham.2019}, \cite{Greuel-review.sings.char.positive} for further references.
   But the $\cA$-case was untouched.

\

Let $\k$ be any field, let  $R_X$ denote the quotient ring of formal power series, $\quots{\k[[x]]}{J}$, resp. analytic,    $\quots{\k\{x\}}{J}$,
 resp. algebraic,$\quots{\k\bl x\br}{J}$,
  see \S\ref{Sec.Notations.Rings.Germs}.
 Accordingly we have the (formal/analytic/algebraic) scheme-germ, $X:=Spec(R_X)$.
  Consider the (formal/analytic/algebraic) map-germs, $\MapX$. The groups  $\cR,\cK,\cA$ act  on $\MapX$. In \cite{Kerner.Group.Orbits}
   I have studied the group orbits, $\cG f$, and their (image) tangent spaces, $T_\cG f$.
    I have obtained various ``linearization" results of type ``$T_\cG f$ vs $\cG f$".

The current paper is the next step. I construct the theory of unfoldings for $\MapX$, the beginning of the theory of stable maps,
 and prove several  local Torelli-type theorems  (the Mather-Yau/Gaffney-Hauser theorems in zero and positive characteristic).

\subsection{The structure/contents of the paper}\label{Sec.Intro.Contents}
  $R_X$ is one of  $\quots{\k[[x]]}{J}$, $\quots{\k\{x\}}{J}$, $\quots{\k\bl x\br}{J}$ and $\cG\in \cR,\cK,\cA$.
\bee[\S 1]\setcounter{enumi}{1}
\item sets the notations for the rings, maps of spaces, groups and their tangent spaces.
 For all the other definitions and results I refer to \cite{Kerner.Group.Orbits}.

Note: through this paper $T_\cG f$ denotes the {\em extended} tangent space. (Classically one writes $T_{\cR^e}f$, $T_{\cK^e}f$, $T_{\cA^e}f$.)
 The classical tangent space is $T_{\cG^{(0)}}f$, for the filtration $\cm^\bullet\cdot \RpX$ on the space of maps.

\item sets  the basic notions of unfoldings (over an arbitrary field): the pullback and $\cG$-equivalence,
 the (infinitesimal) $\cG$-triviality, the (infinitesimal) $\cG$-versality and the (infinitesimal) stability.
  These are copied verbatim from the classical case.

A remark (to avoid any confusion):
\bei
\item In the unfoldings one deforms only the map, not the source (i.e. $X$ remains constant), unlike e.g. \cite{Mond-Montaldi}.
\item Except for \S\ref{Sec.Triviality.of.Unfolding.Nonlocal.Base} we consider the germs of unfoldings, the unfolding base being $(\k^r,o)$.
 Therefore (unlike e.g. \cite{Greuel-Nguyen} or \cite[pg.234]{Gr.Lo.Sh}) we work with versal unfoldings and
  do not introduce complete unfoldings.
\eei

\item treats the triviality of unfoldings.
\bee[\S4.1]
\item gives examples  showing that the classical Thom-Levine criterion cannot hold in positive characteristic. Such pathological unfoldings are
  ``$\cG$-inseparable". The relation to the (in)separability of the group orbit map is explained in \S\ref{Sec.Separable.Unfoldings}.

\item 
  contains the Thom-Levine theorem, ``A separable unfolding is locally trivial iff it is infinitesimally-trivial, i.e. $\di_t f_t\in T_{\cG_t} f_t$".
 (To repeat, in zero characteristic all unfoldings are separable.)

\item
gives an application of this triviality criterion to the  factorization, in $char(\k)=0$:
\bee
\item If the germ $X$ admits $r$ vector fields that are linearly independent at the origin then $X\cong \tX\times(\k^r,o)$.
\item Moreover, suppose a map $f:X\to (\k^p,o)$ ``does not see"  these vector fields, i.e. $T_\cR f=T_{\cR^{(o)}}f$, resp. $T_\cK f=T_{\cK^{(o)}}f$, resp.
 $T_\cA f=T_{\cR^{(o)}}f+T_{\cL}f$. (See \S\ref{Sec.Notations.Groups.Tangent.Spaces}.ii.)
 Then $f$ factorizes as $X\to \tX\to (\k^p,o)$, i.e. is $\cG$-equivalent to the pullback of a map $\tX\to (\k^p,o)$.
\eee
Here part a. is well known for $\C$-analytic germs, but the algebraic version (over $\quots{\k\bl x\br}{J}$) seems new.
 Part b. seems new in all cases.

\item
 is about the rank of a map $f:X\to (\k^p,o)$ and   the corresponding ``preliminary form", splitting $f$ into its linear part and the part of order$\ge2$.

\item
The unfolding trivializations of \S\ref{Sec.Triviality.of.Unfolding.Thom-Levine} are local, they hold over the germ of parameter space $(\k^r_t,o)$.
 In the $\k$-analytic case one gets the trivialization over small balls, $o\in Ball_\ep\sset \k^r$. In the algebraic case ($R_X=\quots{\k\bl x\br}{J}$)
  one gets the trivialization over \'etale neighborhoods of $o\in \k^r$. But in some cases one works over a global base, e.g. $R_{X,t}$ is one of $\quots{\k[t][[x]]}{J}$,
   $ \quots{\k[t]\{x\}}{J}$,  $ \quots{\k[t]\bl x\br}{J}$. Then one wants the global trivialization, over the whole $\k^r$. This is possible if one
    restricts to the filtration-unipotent subgroup $\cG^{(1)}<\cG$.
    In \S\ref{Sec.Triviality.of.Unfolding.Nonlocal.Base} we prove: If $\di_t f_t\in T_{\cG^{(1)}_t}f_t$ then
    the (separable) unfolding $(f_t,t)$ is globally trivial.

The weaker condition  $\di_t f_t\in T_{\cG^{(1)}_t}f_t$ implies the weaker statement: all the fibres are $\cG$-equivalent.

\eee

\item treats the (local) versality of unfoldings.
\bee[\S5.1]
\item  
 gives the pre-normal form: any unfolding $(f_t,t)$
 is formally $\cG$-equivalent to the unfolding $(f_o+\sum a_j(t)v_j,t)$, where $\{v_\bullet\}$ go to the generators of the $\k$-vector space $T^1_\cG f_o$.
\item
  As a simple corollary one gets: an unfolding is versal iff it is infinitesimally versal.
\item
  Recall that the orbit of contact group, $\cK f$, is essentially larger than that of the left-right group, $\cA f$. Hence several
   old questions: ``How does a $\cK$-orbit split into $\cA$-orbits?", ``When do these orbits coincide?",
    ``For which maps is $\cA f\sset \cK f$ an open dense subset?".
     In \S\ref{Sec.Versality.Fibration.of.K.unfoldings.into.A} we address the local version:
      How does a $\cK$-trivial      unfolding split into a family of $\cA$-unfoldings?
\item 
       treats the global version:
       ``How does a $\cK$-orbit split into a family of $\cA$-orbits?".
\eee

\item extends several results of Mather on stable maps to the maps $X\to (\k^p,o)$, in arbitrary characteristic.
\bee[\S6.1]
\item A map is locally stable iff it is infinitesimally stable, i.e. $T_\cA f=\RpX$.
 \item Stable maps are unfoldings of their genotypes, i.e. each stable map is $\cA$-equivalent to $(f_o+\sum u_i v_i,u)$, where
  the elements $\{v_\bullet\}$ generate the vector space $(x)\cdot T^1_\cK f_o$.
\item   Stable maps are determined by their local algebras. Namely, two stable maps are $\cA$-equivalent iff their genotypes are $\cK$-equivalent.
\eee

\item extends the classical results of \cite{Mather-Yau}, \cite{Gaffney-Hauser}, \cite{Scherk}  to the case of $\k$ of zero characteristic, and establishes the
 corresponding versions in positive characteristic. For isolated hypersurface singularities ($p=1$, $R_X=\k[[x]]$) and $\cG\in \cR,\cK$
  this was done in \cite{Greuel-Pham.mather-yau.char.positive}. Our versions ($p\ge1$, $R_X$ one of $\quots{\k[[x]]}{J}$,
   $\quots{\k\{x\}}{J}$, $\quots{\k\bl x\br}{J}$,
    and  $\cG\in \cR,\cK,\cA$) do not assume
   that the singularity is isolated. Moreover, the $char(\k)>0$ version is essentially stronger than that of  \cite{Greuel-Pham.mather-yau.char.positive}.
\eee

\section{Notations and conventions}\label{Sec.Notations.Conventions}
 Below we recall only the main notions. See  \cite[\S2,\S3]{Kerner.Group.Orbits} for the full exposition.

\subsection{Rings and germs}\label{Sec.Notations.Rings.Germs}
\bee[\bf i.]
\item  In this paper $R_X$ is one of the rings $\quots{\k[[x]]}{J}$, $\quots{\k\{x\}}{J}$, $\quots{\k\bl x\br}{J}$.
 Here:
 \bei
 \item  $\k$ is any field, for  $\quots{\k\{x\}}{J}$ we assume $\k$ to be normed and complete \wrt its norm.
\item
  $x=(x_1,\dots,x_n)$ and $J\sseteq (x)^2$.
\eei
   We denote the images of $\{x_i\}$ in $R_X$ by the same letters, this causes no confusion.
   E.g. the maximal ideal $\cm\sset R_X$ is generated by $\{x_i\}$, sometimes we write $\cm=(x)$.

\

 When $J\neq 0$  we always assume  the following $jet_0$-condition.
 Take a derivation $\xi\in  Der_\k(R_X)$, present it as $\sum c_i(x)\di_{x_i}$. Suppose $\xi$ vanishes at the origin to the second order,
  i.e. $c_\bullet(x)\in (x)^2\sset R_X$. If $R_X$ is regular, i.e. $J=0$, then the map $\Phi:x\to x+\xi(x)$
   defines a coordinate change, $\Phi\in  Aut_\k(R_X)$.
In the non-regular case the map $x\to x+\xi(x)$ is not an automorphism of $R_X$.

{\bf$\pmb{jet_0}$  assumption:} any derivation $\xi\in Der^{(1)}_\k(R_X)$ induces an automorphism $x\to x+\xi(x)+h(x)$, where $h(x)$ is a higher order term in the following sense:
  $(h(x))\sseteq (\{c_\bullet(x)\})^2\sset R_X$.

This $jet_0$-assumption  holds if $char(\k)=0$, \cite{BGK.20}. But in positive characteristic this is a non-trivial condition.

\item  The source of a map is the scheme-germ $X:=Spec(R_X)$. The target is the germ $(\k^p,o):=Spec(R_Y)$, where $R_Y$ is   $\k[[y]]$, resp. $\k\{y\}$,
 resp. $\k\bl y\br$.
 Fix some local coordinates in the target $(\k^p,o)$, then the space of maps $\MapX$ is identified with the $R_X$-module $\cm\cdot\RpX$.

Given a map $f:X\to (\k^p,o)$ take its dual $f^\sharp:R_Y\to R_X$ and $f^\sharp:R^{\oplus p}_Y\to \RpX$.
 For a submodule $\La_Y\sset R^{\oplus p}_Y$ we distinguish between the image $f^\sharp(\La_Y)\sset \RpX$ and the pullback
 $f^*(\La_Y):=R_X\cdot f^\sharp(\La_Y)\sset \RpX$.

\item
To work with deformations/unfoldings we take $\K$ as one of the rings $\k[[t]]$,
   $\k\{t\}$,  $\k\bl t\br$, here $t=(t_1,\dots,t_r)$.
    Accordingly we extend the $\k$-algebras $R_X,R_Y$ to
      the $\K$-algebras, $R_{X,t}$ is one of $\quots{\K[[x]]}{J}$, $\quots{\K\{x\}}{J}$, $\quots{\K\bl x\br}{J}$, and $R_{Y,t}$ is one of
 $\K[[y]]$, $\K\{y\}$, $\K\bl y\br$.

\item We often use the Nakayama lemma over a local ring: if $M= C+\cm\cdot M$ (for a finitely generated $M$-module and its subset $C$) then $M= C$.

\item The nested version of Artin approximation, \cite[\S5.2]{Rond}.
\eee

\subsection{The groups $\cG\circlearrowright \MapX$ and tangent spaces $T_\cG$}\label{Sec.Notations.Groups.Tangent.Spaces}
\bee[\bf i.]
\item
The group of right equivalences consists of the coordinate changes  in the source, $\cR:=Aut_X:=Aut_\k(R_X)$.
 The action $\cR\circlearrowright \MapX$ is given by $f\to \Phi_X(f):=f\circ\Phi_X^{-1}$.

 Similarly the group of left equivalences, $\cL:=Aut_{(\k^p,o)}:=Aut_\k(R_Y)\circlearrowright \MapX$, acts by $f\to \Phi_Y(f)$.
 Then group of left-right transformations is $\cA:=\cL\times \cR$.

Recall that the orbits of the contact group,  $\cK:=\cC\rtimes \cR$,
 coincide with the orbits of a much smaller group, $\cK^{lin}:=GL(p,R_X)\rtimes Aut_X$,
 see \cite[\S3.2]{Kerner.Group.Orbits}.

 When working with families, i.e. $R_{X,t}$,  $R_{Y,t}$,  are $\K$-algebras, we take the $\K$-linear automorphisms, $Aut_{X,t}:=Aut_\K(R_{X,t})$,
  $Aut_{Y,t}:=Aut_\K(R_{Y,t})$.

 The tangent spaces to the groups $\cR,\cK,\cA$ are defined (and studied) in \cite[\S3]{Kerner.Group.Orbits}.
   As we mentioned in $\S\ref{Sec.Intro.Contents}$,
   $T_\cG$   denotes the {\em extended}  tangent space:
 \beq
 T_\cR f=Der_X f,\quad\quad\quad  T_\cL f=f^\sharp(\RpY),
 \quad\quad\quad T_\cK f=T_\cR f+(f)\cdot \RpX , \quad\quad\quad T_\cA f=T_\cR f+T_\cL f.
\eeq
\item
Take the filtration $I^\bullet\cdot \RpX$ on $\MapX\cong \cm\cdot \RpX$.  We get the associated filtrations,
 the subgroups  $\cG=\cG^{(-1)}\ge \cG^{(0)}\ge \cG^{(1)}\ge\cdots$  and the submodules
  $T_\cG=T_{\cG^{(-1)}}\supseteq T_{\cG^{(0)}}\supseteq T_{\cG^{(1)}}\supseteq\cdots$.
 Here:
\beq
T_{\cR^{(j)}}:=\{\xi\in Der_X|\ \xi(I^\bullet)\sseteq I^{\bullet+j},\ \forall\ \bullet\ge0\}, \quad\quad
T_{\cL^{(j)}}f= f^\#((y)^{j+1}\cdot \RpY),
\eeq
\[
   T_{\cK^{(j)}}f=T_{\cR^{(j)}}f+I^{j-1}(f)\cdot \RpX,\quad\quad\quad   T_{\cA^{(j)}}f=T_{\cR^{(j)}}f+T_{\cL^{(j)}}f.
\]
E.g. if $R_X$ regular (i.e. $J=0$), and $I=\cm$, then the space $T_{\cG^{(0)}}$ is the ``classical" tangent space.

\item
Take a $\cK$-finite map $f\in \cm\cdot\RpX$. Fix some elements $\{v_\bullet\}$ in $\cm\cdot \RpX$ that go to a basis of the $\k$-vector
 space $\cm\cdot T^1_\cK f$, here $T^1_\cK f:=\quots{\RpX}{T_\cK f}$.
 Thus  $Span_\k\{v_\bullet\}+ T_\cR f+(f)\cdot \RpX+\k^p=\RpX$.
The following presentation is used often.
\bel\label{Thm.TA.vs.TK}
\bee
\item $Span_{R_Y}\{v_\bullet\}+T_\cA f=\RpX$.\quad\quad
\item \text{If } $(f)\sseteq\cm^2\sset R_X$ \text{ then } $Span_{R_Y}\{v_\bullet\}+T_\cR f+  T_{\cL^{(0)}} f =\cm\cdot\RpX$.
\eee\eel
\bpr
{\bf a.} We have $\quots{(Span_{R_Y}\{v_\bullet\}+T_\cA f)}{T_\cR f}+(f)\cdot\quots{\RpX}{T_\cR f}=\quots{\RpX}{T_\cR f}$.
 Consider $\quots{\RpX}{T_\cR f}$ as an $R_Y$-module. It is finitely generated, because  $f$ is $\cK$-finite, \cite[\S3]{Kerner.Group.Orbits}.
   Therefore (by Nakayama over $R_Y$) we get
  $\quots{Span_{R_Y}\{v_\bullet\}+T_\cA f}{T_\cR f} =\quots{\RpX}{T_\cR f}$. Hence $Span_{R_Y}\{v_\bullet\}+T_\cA f=\RpX$.

Now b.   follows from a.
\epr

\eee

\subsection{Changing the base field}
Let $R_X$ be one of $\quots{\k[[x]]}{J},   \quots{\k\{x\}}{J},  \quots{\k\bl x\br}{J}$ and take a(ny) field extension $\k\hookrightarrow\K$.
 Take the ring extension (in the formal case) $R_{X,\K}:=\quots{\K[[x]]}{J}$.

 Given  $\cG\in \cR,\cK,\cA$ we get the group $\cG_\K$,
 e.g. $\cR_\K=Aut_\K(R_{X,\K})$. Take the filtration $M_\bullet:=I^\bullet\cdot \RpX$.

\bel\label{Thm.Equivalence.Change.base.field}  Let $R_X=\quots{\k[[x]]}{J}$ for $\cG=\cA$ or $R_X\in \quots{\k[[x]]}{J},
   \quots{\k\{x\}}{J},  \quots{\k\bl x\br}{J}$ for
 $\cG\in \cR,\cK$.
 Take two maps $f_0,f_1\in \cm\cdot\RpX$.
\bee
\item
Suppose $char(\k)=0$. Let $j\ge1$. If $f\stackrel{\cG^{(j)}_\K}{\sim}\tf$ then $f\stackrel{\cG^{(j)}}{\sim}\tf$.
\item Suppose either the algebraic
  closure $\bar\k$ is uncountable or ($char(\k)=0$ and $\sqrt{I}=\cm$).
Let $j\ge-1$.
 If $f_0\stackrel{\cG^{(j)}_\K}{\sim}f_1$ then $f_0\stackrel{\cG^{(j)}_\bk}{\sim}f_1$.
\eee\eel
\bpr It is enough to consider only the case of   $R_X=\quots{\k[[x]]}{J}$. For the rings $\quots{\k\{x\}}{J}$,  $\quots{\k\bl x\br}{J}$ and the groups  $\cG\in \cR,\cK$ one invokes the Artin approximation,
 \S\ref{Sec.Notations.Rings.Germs}.
\bee
\item  Suppose $g_\K f_0=f_1$ for some $g_\K\in \cG^{(1)}_\K$. Present this group element as the exponential from the tangent space,
 $g_\K=e^{\xi^\K_d}$, here  $\xi^\K_d\in T_{\cG^{(1)}_\K}=\K\hat\otimes T_{\cG^{(1)}}$. (See e.g. \cite[\S3.2]{BGK.20}.)
 One has $ord(\xi^\K_d (f))= ord(f)+d$    for some $d\ge1$.
  Therefore $e^{\xi^\K_d}f-f-\xi^\K_d f\in \K\hat\otimes M_{ord(f)+d+1} $.
   As $e^{\xi^\K_d}f,f\in \RpX$, we get
   $\xi^\K_d f\in \RpX+\K\hat\otimes M_{ord(f)+d+1} $.
    Pass to the quotient vector space,
   \beq
   [\xi^\K_d f]\in \quot{\RpX+\K\hat\otimes M_{ord(f)+d+1} }{\K\hat\otimes M_{ord(f)+d+1} }\sset
   \K\hat\otimes\quot{\RpX }{ M_{ord(f)+d+1}}.
\eeq
Recall the general fact: for a vector subspace $V_\k\sset W_\k$ one has $(\K\hat\otimes V_\k)\cap W_\k=V_\k$. Therefore
 $\xi^\K_d f\in T_{\cG^{(1)}}f+\K\hat\otimes M_{ord(f)+d+1}$. Thus we can expand
 $\xi^\K_d=\xi_d+\xi^\K_{d+1}$,\  where $\xi_d\in T_{\cG^{(1)}}$, $\xi^\K_{d+1}\in \K\hat\otimes T_{\cG^{(1)}}$,   and
  $ord(\xi_d f)\ge   ord( f)+d$,    $ord(\xi^\K_{d+1} f)\ge   ord(f)+d+1$.
 Replace $\tf$ by $e^{-\xi_d}\tf$, and iterate.

Eventually we get the infinite product $g:=\lim_d (e^{\xi _d}\cdots e^{\xi _1})\in \cG^{(1)}$.
 We claim: this limit exists. Indeed, $\xi_d f\in T_{\cG^{(1)}}f\cap I^d\cdot\RpX$. Then by the Artin-Rees type property, \cite[\S6]{Kerner.Group.Orbits},
  we get:  $\xi_d f\in T_{\cG^{(j_d)}}f$, with $\lim_d j_d=\infty$.

Altogether we have $\tf=g\cdot f$. This proves the statement for $R_X=\quots{\k[[x]]}{J}$.

\item (We show only the case of $\cA$-equivalence, the other cases are similar.) We should resolve the condition $\Phi_Y\circ f_0\circ\Phi_X=f_1$.
 Take the Taylor expansions, $\Phi_X(x)=\sum_{|m|\ge1} C^{(x)}_m x^m$ and  $\Phi_Y(y)=\sum_{|m|\ge1} C^{(y)}_m y^m$. Here $m$ is the multi-index, and
  $\{C^{(x)}_m\}$,   $\{C^{(y)}_m\}$ are unknowns. By comparing the coefficients of the monomials, the condition
    $\Phi_Y\circ f_0\circ\Phi_X=f_1$ (and $\Phi_X(J)=J\sset \k[[x]]$) is transformed into a countable system of polynomials equations, $\{P_j(\{C^{(x)}_*\},\{C^{(y)}_*\})=0\}_j$.
     Each polynomial here is in a finite number of variables.

     This system is solvable over $\K$. We claim: each finite subsystem of $\{P_j(\{C^{(x)}_*\},\{C^{(y)}_*\})=0\}_j$ is solvable over $\bar\k$.
 Indeed, this subsystem       defines a subscheme in an affine space over $\bk$.
       If this finite subsystem is not solvable then this subscheme has no points over $\bk$.
       By Hilbert Nullstellensatz we get: the defining ideal of this subscheme is the whole polynomial ring of the affine space.
        And then this finite subsystem can have
        no solutions over $\K$.

        Finally we apply the assumptions on $\k$.
\bei
\item
Suppose $\bar\k$ is uncountable. In this case a countable system of polynomial equations over $\bk$ is solvable iff each finite subsystem if solvable,
 \cite[Theorem 5]{Popescu-Rond}.  Thus the solvability over $\K$ implies that over $\bar\k$.
\item If $char(\k)=0$ and $\sqrt{I}=\cm$, then it is enough to take only finite expansions in $\Phi_X(x)=\sum_{|m|\ge1} C^{(x)}_m x^m$,
   $\Phi_Y(y)=\sum_{|m|\ge1} C^{(y)}_m y^m$. Resolving the corresponding equations one gets $f\stackrel{\cG_{\bar\k}}{\sim}\tf$
    with $\tf\stackrel{\cG^{(1)}_\K}{\sim}f$. Now apply part 1. \epr
\eei
\eee

\beR
The first statement fails in positive characteristic.  For example, let $f(x)=x^p\in \k[[x]]$, where $char(\k)=p$.
Then $x^p\stackrel{\bar\k\otimes\cR}{\sim}x^p+ax^{2p}$  for any $a\in \k$. On the other hand,
 suppose the Frobenius morphism is non-surjective (i.e. the field $\k$
 is non-perfect). Then $x^p\stackrel{\cR}{\not\sim}x^p+ax^{2p}$ for any $a\in \k\smin \k^p$.
\eeR

\section{The basic notions of unfolding}\label{Sec.Unfoldings.Definitions}
Let $R_X$ be one of $\quots{\k[[x]]}{J}$, $\quots{\k\{x\}}{J}$, $\quots{\k\bl x\br}{J}$, accordingly take $\K$ and
 $R_{X,t}$, see \S\ref{Sec.Notations.Rings.Germs}.
 \subsection{}
\bed
An unfolding of a map $f:X {\to} (\k^p,o) $ is the map $F:X\times(\k^r,o){\to}(\k^{p},o)\times (\k^{r},o)$ of the form
 $F(x,t)=(f_t(x),t)$,
 i.e.  an element $(f_t,t)\in (\cm+(t))\cdot R^{\oplus (p+r)}_{X,t}$, satisfying $f_{o}=f$.
\eed

   Take the actions $ \cR,\cL \circlearrowright \cm\cdot \RpX$.
    Accordingly we have the actions $ \cR_t,\cL_t \circlearrowright (\cm+(t))\cdot \RpXt$, where
\beq
\cR_t:=\{\Phi_{X,t}|\ \Phi_{X,o}=Id_X\}< Aut_\K(R_{X,t}),\quad\quad\quad\quad
 \cL_t:=\{\Phi_{Y,t}|\ \Phi_{Y,o}=Id_Y\}< Aut_\K(R_{Y,t}).
\eeq
Geometrically we have the coordinate changes $\cR_t\circlearrowright X\times(\k^r,o)$ and  $\cL_t\circlearrowright (\k^{p+r},o)$
that restrict to identities on  the central fibres $X\times \{o\}$, $\k^p\times \{o\}$, and preserve  all the $t=const$ slices,
  $X\times\{t\}$, $(\k^p,o)\times\{t\}$.
 Equivalently:  the  elements of $\cR_t$, $\cL_t$ are unfoldings of identity maps, $(x,t)\to (x,t)$ and $(y,t)\to (y,t)$.

Similarly one extends the groups $\cA,\cK$, $\cK^{lin}=GL(p,R_X)\rtimes \cR$ to:
\beq
\cA_t:=\cL_t\times \cR_t\circlearrowright (x,t)\cdot\RpXt, \quad    \cK_t:=\cC_t\rtimes \cR_t\circlearrowright (x,t)\cdot\RpXt,   \quad
\cK^{lin}_t:= GL(p,R_{X,t})\rtimes \cR_t\circlearrowright (x,t)\cdot\RpXt.
\eeq
These groups act on families of maps, $\cG_t\circlearrowright Maps(X\times(\k^r,o),(\k^{p+r},o))$. Here the groups $\cC_t$, $GL(p,R_{X,t})$
  preserve the origin of $(\k^p,o)$ (for each $t$). The action $\cR_t$ (resp. $\cL_t$) does not preserve the origin of $(\k^n,o)$ (resp. $(\k^p,o)$),
 i.e. $\Phi_{X,t}(x)\not\sseteq (x)$ and $\Phi_{Y,t}(y)\not\sseteq (y)$.
  To preserve the origin(s) one considers the subgroups  $\cG^{(o)}_t$ and subspaces
 $T_{\cG^{(o)}_t}$ for the filtration $(x)^\bullet$, see \ref{Sec.Notations.Groups.Tangent.Spaces}.

  The actions of $\cR_t$, $\cK^{lin}_t\circlearrowright (x,t)\cdot\RpXt$ are $t$-linear. The actions of  $\cK_t$, $\cL_t$, $\cA_t$ are not $t$-linear
   (neither additive nor $\k$-multiplicative).

For  $R_X=\k\{x\}$ (with $\k\in \R,\C$) or $C^\infty(\R^n,o)$, we get the classical unfolding notions.
 The map is $(\k^n,o)\stackrel{f}{\to}(\k^p,o)$ and its $t$-unfolding is
 $(\k^n\times\k^r,o)\stackrel{F}{\to}(\k^p\times\k^r,o)$, $(x,t)\to (f_t(x),t)$.

\

The (extended) tangent spaces to these groups are:
\beq
T_{\cR_t}:=Der_\K(R_{X,t}),\quad  T_{\cL_t}:=Der_\K(R_{Y,t}), \quad     T_{\cA_t}:=T_{\cL_t}\oplus T_{\cR_t}, \quad
 T_{\cK_t}:=T_{\cR_t}\oplus Mat_{p\times p}(R_{X,t}).
\eeq
Accordingly we have the image tangent spaces $T_{\cG_t}f_t\sseteq \RpXt$ and $T_{\cG_t}F\sseteq R^{p+r}_{X,t}$.

The standard notions of unfoldings are extended (verbatim) from the classical case, \cite{AGV},
 \cite[III.2]{AGLV}, \cite[\S4,\S5]{Mond-Nuno}, \cite[Chapter XIV]{Martinet.1982}.
 Let $\cG$ be one of $\cR,\cK,\cA$, accordingly $\cG_t\in \cR_t,\cK_t,\cA_t$.
\bed\bee[\bf 1.]
\item Two unfoldings $F, \tF\in (\cm+(t))\cdot R^{\oplus (p+r)}_{X,t}$ of $f\in \cm\cdot \RpX$ are called $\cG_t$-equivalent if $g_t\cdot f_t=\tf_t$ for some $g_t\in \cG_t$.

\item\bei
\item
An unfolding $F$ is called $\cG$-trivial if it is $\cG_t$-equivalent to the constant unfolding, i.e.  $F\sim (f_{o},t)$.

\item
An unfolding $F$ is called infinitesimally-$\cG$-trivial if $Span_\k(\di_{t_1}f_t,\dots, \di_{t_r}f_t)\sseteq T_{\cG_t}f_t$.
\eei

\item
\bei
\item A map $f\in Maps (X,(\k^p,o)  )$ is called $\cA$-stable if all its unfoldings are $\cA$-trivial.
\item A map $f\in Maps (X,(\k^p,o)  )$ is called infinitesimally-$\cA$-stable if $T_{\cA}f=\RpX$.
\eei

\item The pull-back of an unfolding $(f_t(x),t)\in (\cm+(t))\cdot R^{\oplus(p+r)}_{X,t}$ via the map $\phi:\ (\k^\tr_{\tt},o)\to (\k^r_t,o)$ is the unfolding
  $(f_{t(\tt)}(x),\tt)\in (\cm+(\tt))\cdot R^{\oplus(p+\tr)}_{X,\tt}$.
   Here $t(\tt)$ is the image   $\phi^\sharp(t)$ for the map $\phi^\sharp:\ \K^r_t\to \K^\tr_\tt$.
\item
\bei
\item  An unfolding $F\in (\cm+(t))\cdot R^{\oplus (p+r)}_{X,t}$ of $f$ is called $\cG$-versal
 if any other unfolding of $f$ is $\cG_t$-equivalent to a pullback of $F$.
\item
An unfolding $F\in (\cm+(t))\cdot R^{\oplus (p+r)}_{X,t}$ of $f$ is called infinitesimally-$\cG$-versal if
 \[Span_\k(\di_{t_1}f_t|_{t=o},\dots, \di_{t_r}f_t|_{t=o})+ T_{\cG}f_o=\RpX.\]
 Namely, the elements $\di_{t_1}f_t|_{t=o},\dots, \di_{t_r}f_t|_{t=o}$ are sent to generators of the vector space $T^1_\cG f:=\quots{\RpX}{T_\cG f_o}$.
\eei
\eee
\eed

\beR
\bee[\bf i.]
 \item
The stability notions (in part 3) are introduced  for $\cA$-equivalence only. They are not useful for  $\cR$,$\cK$ equivalences.
  Indeed, any   infinitesimally $\cR$ or $\cK$-stable
 map is necessarily equivalent to $f(x)=(x_1,\dots,x_p)$. (Thus in particular $p\le n$.)

\item These definitions are over local  rings. Geometrically one works in the infinitesimal neighborhood of $o\in X$. For the ring $\quots{\k\{x\}}{J}$
 one can pass (in the standard way) to small neighborhoods. E.g., suppose the unfolding $F$ is trivial, i.e. $F\sim(f_o,t)$. Then there exist
  small balls, $Ball_\ep\sset X$ and $o\in Ball_\de\sset \k^r$, with $0<\de\ll\ep$, and an analytic family of elements $\{g_t\in \cG_t\}_{t\in Ball_\de}$,
    satisfying (in $Ball_\ep$): $g_t f_t=\tf_t$ for each $t\in Ball_\de$.

For the ring $\quots{\k\bl x\br}{J}$ the small neighborhoods are the \'etale covers.
  If $F$ is $\cG$-trivial then there exists an \'etale map $\cU_\tt\stackrel{\phi}{\to}\k^r_t$,
 with $\phi(\cU)$ a Zariski open neighborhood of $o\in \k^r$, such that the pullback $\phi^\sharp(F)=(f_{t(\tt)},\tt)$ is (globally) trivial over $\cU$.

\item Take two maps $f_o,\!\tf_o\!\in\!\MapX$\! and their unfoldings $F,\!\tF$.\! If $F\!\!\stackrel{\cG}{\sim}\!\!\tF$\! (as unfoldings)
 then $f_o\!\!\stackrel{\cG}{\sim}\!\!\tf_o.$
 \eee
\eeR

\subsection{}
\bel\label{Thm.Unfolding.Properties.Invariant.under.G.action}
The conditions of (infinitesimal)  $\cG$-triviality, (infinitesimal) $\cA$-stability, (infinitesimal) $\cG$-versality are
 preserved under the $\cG$-equivalence of unfoldings.
\eel
Namely, if $F\stackrel{\cG_t}{\sim}\tF$, and $F$ has one of these properties, then so does $\tF$.
\\\bpr  The invariance of the triviality/stability/versality under $\cG$-equivalence is tautological.
\bei
\item ({\em infinitesimal  $\cG$-triviality})
 We prove: the condition $\di_t f_t\in T_{\cG_t}f_t$ is preserved under   $\cG_t$-transformations, i.e.
 $\di_t (g\cdot f_t)\in T_{\cG_t}(g\cdot f_t)$ for any $g\in \cG_t$.

 Suppose $J=0$.
 Below $f$ is a power series in $x$, while $f'$ denotes the $x$-derivative.
\bei
\item{\bf  The $\cR$-case.} We should prove: $\di_t (\Phi_{X,t}( f_t))\in T_{\cR_t}(\Phi_{X,t}( f_t))$.
 First we verify: $\Phi_{X,t}(T_{\cR_t}f_t)=T_{\cR_t}(\Phi_{X,t}(f_t) )$. This is the chain rule:
 \beq\label{Eq.Invariance.R.infinites.trivial}
 T_{\cR_t}(\Phi_{X,t}(f_t))\ni \xi_X(f_t\circ\Phi_{X,t})=\sum\nolimits_j \di_j f_t|_{\Phi_{X,t}(x)}\cdot \xi_X(\Phi_{X,t}(x)_j)\in
\Phi_{X,t}( T_{\cR_t}f_t).
\eeq
Now we observe: $\di_t (f_t\circ\Phi_{X,t})=\di_t f_t|_{ \Phi_{X,t}}+f'_t|_{ \Phi_{X,t}}\cdot \di_t  \Phi_{X,t}\in \Phi_{X,t}(T_{\cR_t}f_t)+T_{\cR_t}\Phi_{X,t}(f_t) $. Hence the statement.
\item{\bf  The $\cK$-case.} Assuming $\di_t f_t\in T_{\cK_t}f_t$ we have:
\beq\label{Eq.Invariance.K.infinites.trivial}
\di_t (U\cdot (f_t\circ \Phi_{X,t}))=(\di_t U)\cdot (f_t\circ \Phi_{X,t})+U\cdot  \di_t   f_t|_{\Phi_{X,t}}+
  U\cdot f'_t|_{\Phi_{X,t}}\cdot \di_t \Phi_{X,t} \in T_{\cK_t}(U\cdot (f_t\circ\Phi_{X,t})).
\eeq
\item{\bf  The $\cA$-case.}  Assuming $\di_t f_t\in T_{\cA_t}f_t$ we have:
\beq\label{Eq.Invariance.A.infinites.trivial}
\di_t(\Phi_{Y,t}\circ f_t\circ\Phi_{X,t})=(\di_t \Phi_{Y,t})|_{f_t\circ\Phi_{X,t}}+\Phi'_{Y,t}|_{f_t\circ\Phi_{X,t}}\cdot \di_t f_t|_{\Phi_{X,t}}
+\Phi'_{Y,t}|_{f_t\circ\Phi_{X,t}}\cdot f'_t|_{\Phi_{X,t}}\cdot \di_t \Phi_{X,t}\in
\eeq
\[
\in T_{\cL_t} (\Phi_{Y,t}\circ f_t\circ\Phi_{X,t})+\Phi'_{Y,t}|_{f_t\circ \Phi_{X,t}}\cdot T_{\cA_t}f|_{\Phi_{X,t}}+\Phi'_{Y,t}|_{f_t\circ \Phi_{X,t}}
 \cdot T_{\cR_t} f_t|_{\Phi_{X,t}} \sseteq  T_{\cA_t} (\Phi_{Y,t}\circ f_t\circ\Phi_{X,t}).
\]
\eei

For $J\neq 0$ one repeats these arguments for representatives of $f_t,\Phi_{X,t}$.
 Namely, let $S$ be one of $\k[[x]]$, $\k\bl x\br$, $\k\{x\}$ and take a representative $\tf_t\in (x,t)\cdot S^{\oplus p}_{x,t}$ of $f_t\in \RpXt$.
  The group $Aut_\K(R_t)$ lifts to $Aut_{\K,J}(S_{x,t})$, these are automorphisms that preserve the ideal $J$.
   The tangent space is now $T_\cR:=Der_S(log(J))$. Then proceed as before.

\item ({\em infinitesimal  $\cA$-stability}) $T_\cA(\Phi_Y\circ f\circ\Phi_X)=\Phi_Y'|_{f\circ\Phi_X}\cdot T_\cR f|_{\Phi_X}+T_\cL(f\circ\Phi_X)=
 \Phi_Y'|_{f\circ\Phi_X}\cdot T_\cA f|_{\Phi_X}=\RpX$. Here we use the invertibility of $ \Phi_Y'$.
\item
({\em infinitesimal  $\cG$-versality}) The\! verification\! is\! again\! the chain rule, as in the equations
\eqref{Eq.Invariance.R.infinites.trivial},\eqref{Eq.Invariance.K.infinites.trivial},\eqref{Eq.Invariance.A.infinites.trivial}.
\epr\eei

\

Take a map $f_o\in (x)\cdot\RpX$ and its deformation $f_t\in (x,t)\cdot \RpXt$. Let $\cG\in \cR,\cK,\cA$, accordingly $\cG_t\in \cR_t,\cK_t,\cA_t$.
 Let $\La\sset\RpX$ be a $\k$ or $R_Y$-submodule, finitely generated by $\{v_\bullet\}$. Associate to $\La$ the submodule $\La_t:=\K\{v_\bullet\}$, resp.
  $R_{Y,t}\{v_\bullet\}$.
\bel\label{Thm.Unfoldings.Algebraic.Lemma} (The ``algebraic lemma of unfolding", e.g. \cite[pg. 193]{Martinet.1982})

Suppose $T_\cG f_o+\Lambda=\RpX$, where $\La\sset \RpX$ is a finite-dimensional $\k$-vector subspace or a finitely-generated $R_Y$-submodule.
 Then $T_{\cG_t} f_t+ \Lambda_t=\RpXt$.
\eel
\bpr We have the obvious presentation $\RpXt=T_{\cG_t} f_t+(t)\cdot \RpXt+\Lambda_t$. Thus
 $\quots{\RpXt}{T_{\cG_t} f_t}=(t)\cdot \quots{\RpXt}{T_{\cG_t} f_t}+\quots{\Lambda_t+T_{\cG_t} f_t}{T_{\cG_t} f_t}$.
   Note: $\quots{\K}{(t) }\otimes\quots{\RpXt}{T_{\cG_t} f_t}\cong \quots{\RpX}{T_{\cG} f_o}$, the later quotient being finitely generated over $\k$ (resp. $R_Y$).
 Therefore  $\quots{\RpXt}{T_{\cG_t} f_t}$ is finitely generated by Weierstra\ss\ finiteness, \cite[\S2]{Kerner.Group.Orbits}.
  Finally, Nakayama over $\K$ (resp. $R_{Y,t}$) gives  $\quots{\RpXt}{T_{\cG_t} f_t}= \quots{\Lambda_t+T_{\cG_t} f_t}{T_{\cG_t} f_t}$.
\epr

\section{Triviality of unfolding}\label{Sec.Triviality.of.Unfolding}

\subsection{}\label{Sec.Triviality.of.Unfolding.Bad.Examples}
 As in the classical case,   $R_X=\C\{x\}$, one wants to establish Thom-Levine's criterion, ``triviality
 vs   infinitesimal triviality". This does not hold in   full generality because of the positive characteristic obstructions.
 \bex\label{Ex.Unfoldings.pathology.positive.char} Let $\k$ be a field of characteristic $p$ and take $R_X=\k[[x]]$, $\cG=\cR$.
 \bee[\bf i.]

 \item The unfolding $f_t(x)=x^n+t^p x$ is non-trivial, even though $\di_t f_t=0\in T_{\cR_t} f_t$. Note that $f_t$ is obtained by the $t\to t^p$ base-change
  from $x^n+tx$.

The unfolding $f_t(x)=x^n+t^px^p+tx^d\in \k[[t,x]]$, for $p<n<d$, is non-trivial, with $\di_t f_t\in T_{\cR_t} f_t$.
 And this is not a $t\to t^p$ base change of another unfolding.  Here one can take $d\gg n$, i.e. to ensure $\di_t F\in T_{\cR^{(j)}}F$ for $j\gg1$.

\item Consider the unfolding $f_t(x)=x^p+x^{p+d}(1+t)\in \k[[t,x]]$, with $gcd(d,p)=1$.
  Then $\di_t f_t\in x\cdot T_{\cR_t} f_t$, but $f_t$ is not an $\cR$-trivial unfolding. Indeed, for a coordinate change $x\to x+x\cdot h(t,x)$ one has
  $f_t( x+x\cdot h(t,x))=x^p+x^p\cdot h(t,x)^p+x^{p+d}(1+h(t,x))^{p+d})(1+t)$. If this coordinate change trivializes $f_t$ then we must have
   $ord_x(x^p\cdot h(t,x)^p)>p$. But then necessarily $ord_x(x^p\cdot h(t,x)^p)=ord(x^{p+d}\cdot t)=p+d$, contradicting the divisibility
   $p| ord_x(x^p\cdot h(t,x)^p)$.
\eee
 \eex

\subsection{Infinitesimal vs local triviality of unfoldings}\label{Sec.Triviality.of.Unfolding.Thom-Levine}

In view of example \ref{Ex.Unfoldings.pathology.positive.char} we introduce:
\bed
An unfolding of the map $f_o\in \cm\cdot \RpX$ is called $\cG$-inseparable if $f_t(x)\stackrel{\cG_t}{\sim}f_o(x)+t^{d}\cdot f_{ d}(x)+(t^{d+1})\cdot \RpXt$,
  where $char(\k)\mid d$ and $f_{ d}\not\in T_{\cG_t}f_o$.
\eed
Thus in zero characteristic all unfoldings are separable.
 The name ``separable" is due to the relation to separability of group orbit map, see \S\ref{Sec.Separable.Unfoldings}.

\bthe\label{Thm.Unfolding.Triviality.Local} Let $R_X$ be one of $\quots{\k[[x]]}{J}$, $\quots{\k\{x\}}{J}$, $\quots{\k\bl x\br}{J}$,
 see \S\ref{Sec.Notations.Rings.Germs}.i.
 Let $\cG\in \cR,\cK,\cA$. Take an unfolding $F(x,t)=(f_t(x),t)$, for $t=(t_1,\dots,t_r)\in (\k^r,o)$.
\bee[\bf 1.]
\item
 If $F$ is trivial then it is infinitesimally  trivial, i.e. $\di_{t_1} f_t,\dots, \di_{t_r} f_t\in T_{\cG_t}f_t$.
\item
Suppose $F$ is   infinitesimally-$\cG$-trivial. For $char(\k)>0$ assume that $F$ is   $\cG$-separable.  Then:
\bei
\item
 (for $\cG=\cR,\cK$)  $F$ is trivial.
\item
 (for $\cG=\cA$)  $F$ is formally-$\cA$-trivial. Moreover:
 \bei
 \item if $char(\k)=0$ and $R_X=\quots{\k\{x\}}{J}$, then the unfolding $F$ is $\cA$-trivial.
 \item if  the map $f_o\in \RpX$  is $\cA_t$-finitely determined (as an element of $\RpXt$),   then the unfolding $F$ is $\cA$-trivial.
  \eei
 \eei
\eee
\ethe
\bpr
First we reduce the case of $r$-parameters to the one-parameter unfolding.
\bei
\item (Part 1.) Consider $(f_t,t_r)$ as a one-parameter unfolding of the map $f_t|_{t_r=0}\in R^{\oplus p}_{X,t_1,\dots,t_{r-1}}$.
 Here $R^{\oplus p}_{X,t_1,\dots,t_{r-1}}$ is the algebra over $\K=\k[[t_1\dots t_{r-1}]]$, resp. $\k\{t_1\dots t_{r-1}\}$,
  resp $\k\bl t_1\dots t_{r-1}\br$.
 This unfolding is still   trivial,
 therefore (assuming the case $r=1$) it is infinitesimally trivial, $\di_{t_r}f_t\in T_{\cG_t} f_t$. Repeating this for all $t_i$ we get the statement.
\item (Part 2.) Consider$(f_t,t_r)$ as a one-parameter unfolding of the map  $f_t |_{t_r=0}\in R^{\oplus p}_{X,t_1,\dots,t_{r-1}}$.
  This unfolding is still infinitesimally trivial,
 $\di_{t_r}f_t\in T_{\cG_t}f_t$.
 Therefore (assuming the case $r=1$) it can be trivialized. Namely, $g^{(r)}_t f_t(x)$
  depends on the variables $(t_1,\dots,t_{r-1},x)$ only, for an element $g^{(r)}_t\in \cG_t$ that acts as identity on $t_1,\dots,t_{r-1}$.
 Iterate this argument to get the full trivialization: $(g^{(1)}_t\cdot \cdots\cdot g^{(r)}_t) f_t(x)=f_o(x)$.
\eei

Therefore, below we consider only one-parameter unfoldings.

 The proof of part 1 is characteristic-free.
 For part 2 we give an additional proof in  $char(\k)=0$ case.

\bee[\bf 1.]
\item
\bei
\item {\bf  The $\cR$-case.} Take the trivialization:  $f_t=\Phi_{X,t}(f_o)$ for $\Phi_{X,t}\in Aut_\K(R_{X,t})$
 with $\Phi_{X,o}=Id$,  see \S\ref{Sec.Notations.Rings.Germs}.iii.
 We should prove: $\di_t f_t\in T_{\cR_t}f_t$.

 Define the operator $\xi_X:=\di_t-\Phi_{X,t}\circ\di_t \circ\Phi^{-1}_{X,t}$.
  We claim: this is a $t$-linear derivation, $\xi_X\in Der_\K(R_{X,t})$. First observe that the action $\xi_X\circlearrowright R_{X,t}$ is well defined.
   Moreover, this action is $\k$-linear. As $\Phi_{X,t}$ is $\K$-linear, we get  $\xi_X(\K)=0$. Thus $\xi_X\in End_\K(R_{X,t})$.
     Now we verify the Leibniz rule. For any $a,b\in R_{X,t}$ one has
    $\xi_X(ab)=\di_t(ab)-\Phi_{X,t}\circ\di_t \circ\Phi^{-1}_{X,t}(ab)=\xi_X(a)b+a\xi_X(b)$. Thus $\xi_X\in Der_\K(R_{X,t})=T_{\cR_t}$.

Finally we observe: $T_{\cR_t}f_t\ni \xi_X (f_t)=\di_t f_t-\Phi_t \di_t f_o=\di_t f_t$.

\item{\bf  The $\cK$-case.} (As always, we replace $\cK$ by $\cK^{lin}$.)
 Suppose $f_t=(U\cdot f_o)\circ \Phi_{X,t}$ for some $\Phi_{X,t}\in Aut_\K(R_{X,t})$  with $\Phi_{X,o}=Id$,
 and $U\in GL(p,R_{X,t})$  with $U|_o=\one$.  As in the $\cR$-case we define the operator
  $\xi_X:= \di_t-\Phi_{X,t}\circ\di_t \circ\Phi^{-1}_{X,t}$. As in the $\cR$-case this is a derivation, $\xi_X\in Der_\K(R_{X,t})$.
  Now we verify:
  \beq
  \xi_X (f_t)=\di_t f_t-\Phi_{X,t}\circ\di_t (U\cdot f_o)=\di_t f_t-\Phi_{X,t}(\di_t U\cdot U^{-1})\cdot f_t.
  \eeq
Therefore $\di_t f_t\in T_{\cK_t} f_t$.

\item{\bf  The $\cA$-case.} Present the trivial unfolding in the form $f_t=(\Phi^{-1}_{Y,t},\Phi^{-1}_{X,t})(f_o)$.
 Define the operator $\Psi\circlearrowright (\cm+(t))\cdot \RpXt$ as follows:
 \beq
 \Psi(h):=(\Phi'_{Y,t})^{-1}|_{h}\cdot  \Phi^{-1}_{X,t}\circ \frac{\di}{\di t}\circ (\Phi_{Y,t},\Phi_{X,t})(h).
 \eeq
  Here $ \Phi'_{Y,t} $ is the derivative operator and  $(\Phi'_{Y,t})^{-1}$ is its inverse.
We observe: $\Psi(f_t)=0\in \RpXt$. Now we check the $\Psi$-action on $\K$. Take $h(t)\in \K^{\oplus p}$ then:
\beq
\Psi(h(t))=(\Phi'_{Y,t})^{-1}|_{h(t)}\cdot     \frac{\di}{\di t}   \Phi_{Y,t} (h(t))=
(\Phi'_{Y,t})^{-1}|_{h(t)}\cdot\Big[(\di_t  \Phi_{Y,t} )|_{h(t)}+\Phi_{Y,t}'|_{h(t)}\cdot \frac{\di h(t)}{\di t}\Big]=
\eeq
\[
=\xi_Y(y)|_{h(t)}+\frac{\di h(t)}{\di t}\in R^{\oplus p}_{Y,t}.
\]
Here $\xi_Y:=(\Phi'_{Y,t})^{-1} \cdot (\di_t  \Phi_{Y,t} )\in T_{\cL_t}$. Therefore $(\Psi-\xi_Y-\frac{d}{dt})(\K^{\oplus p})=0$.
Moreover, for $h\in (\cm+(t))\RpXt$ one has:
\beq
(\Psi-\xi_Y-\frac{\di}{\di t})h|_{(x,t)}=(\Phi'_{Y,t})^{-1}|_{h}\cdot  \Phi^{-1}_{X,t}\circ \frac{\di}{\di t} (\Phi_{Y,t}|_{\Phi_{X,t}(h)})-\xi_Y(h )-
  \frac{\di}{\di t} h =
\eeq
\[
=\Phi_{X,t}^{-1}\circ \frac{\di}{\di t}\Phi_{X,t}|_h-\frac{\di }{\di t}h =h'|_{(x,t)}\cdot (\Phi_{X,t}^{-1}\circ\frac{\di }{\di t}\Phi_{X,t})\in T_{\cR_t} h.
\]
Thus we get the derivation  $\xi_X:=\Psi-\xi_Y-\frac{\di }{\di t}\in T_{\cR_t}$.
 Finally, $0=\Psi(f_t)=(\xi_Y+\xi_X+\di_t )f_t$ gives:   $\di_t f_t=-(\xi_Y+\xi_X)f_t\in T_{\cA_t}f_t$.
\eei

\item {\bf The case $\mathbf{char(\k)=0}$, $\cG\!\in\! \cR,\cK$.}
First we establish the statement in the formal case, $R_{X,t}\!=\!\quots{\k[[x,t]]}{J}$.
\bei
\item {\bf  The $\cR$-case.} Suppose $\di_t f_t=\xi_X f_t$, for a derivation $\xi_X\in Der_{\k[[t]]}(R_{X,t})$.
Then $(\di_t-\xi_X)^j f_t=0$ for each $j\ge1$.

Extend the ring $R_{X,t}$ by  a new formal variable  $\tt$ to $R_{X,t,\tt}=\quots{\k[[x,t,\tt]]}{J}$.
 In this ring we have:   $e^{\tt\cdot (\di_t-\xi_X)}f_t=f_t$. (Note that the derivation $\tt\cdot (\di_t-\xi_X)$ is
  filtration-nilpotent, see \cite[\S2.1]{Kerner.Group.Orbits})
 Therefore  one has the identity:
 \beq\label{Eq.Proof.Trivial.vs.Infinit.Triv.R.case}
 e^{-\tt\cdot \di_t}\cdot e^{\tt\cdot (\di_t-\xi_X)}f_t= e^{-\tt\cdot \di_t} f_t=f_{t-\tt}.
\eeq
Now we use the Baker-Campbell-Hausdorff formula, \S 2.1.ix  of \cite{Kerner.Group.Orbits}:
 $ e^{-\tt\cdot \di_t}\cdot e^{\tt\cdot (\di_t-\xi_X)}=e^{\sum \tt^l\cdot p_l(  \di_t,\di_t-\xi_X)}$.
 Note that the commutator of derivations is a derivation (of first order).

Therefore  $ e^{-\tt\cdot \di_t}\cdot e^{\tt\cdot (\di_t-\xi_X)}=e^{\tt\cdot \zeta}$, for a derivation $\zeta\in  Der_{\k[[t,\tt]]}(R_{X,t,\tt})$.
 In particular, $\tt\cdot\zeta\in \tt\cdot T_{\cR_{t,\tt}}$.
   Therefore we get the coordinate change  $e^{\tt\cdot\zeta}\in \cR_{X,t,\tt}$, which satisfies:
     $e^{\tt\cdot\zeta} f_t=f_{t-\tt}$.
 This equality is an identity in $R_{X,t,\tt}$.
     We can restrict this identity to $\tt=t$, as $\zeta$ is $\k[[t,\tt]]$-linear. We get the trivialization: $e^{ t \zeta|_{\tt=t}}f_t=f_o$.
  Note that $e^{ t\cdot\zeta|_{\tt=t}}\in \cR_t$, i.e. is the unfolding of identity  $Id\in \cR$.

\item{\bf  The $\cK$-case.} Suppose $\di_t f_t=\xi_X (f_t)+u\cdot f_t\in T_{\cK_t}f_t$, here $\xi_X\in Der_{\k[[t]]}(R_{X,t})$,
 $u\in Mat_{p\times p}(R_{X,t})$. As in the $\cR$-case we get    $e^{\tt\cdot (\di_t-\xi_X-u)}f_t=f_t$ in the ring $R_{X,t,\tt}$. Then, as before,
  we have $e^{-\tt\di_t}\cdot e^{\tt\cdot (\di_t-\xi_X-u)}f_t=f_{t-\tt}$. Using the BCH-formula, as before,  we get the identity:
\beq\label{Eq.Proof.Trivial.vs.Infinit.Triv.K.case}
 e^{\tt\cdot (\zeta_X-\cu)}f_t=f_{t-\tt},   \text{ where }   \zeta\in Der_{\k[[t,\tt]]}(R_{X,t,\tt})  \text{ and }
  \cu\in  Mat_{p\times p}(R_{X,t,\tt}).
  \eeq
 In particular $\tt\cdot (\zeta_X-\cu)\in \tt\cdot T_{\cK_{t,\tt}}$.
  As before, one substitutes $\tt=t$ to get the trivialization:   $ e^{t\cdot (\zeta_X-\cu)|_{\tt=t}}f_t=f_o$.
    Here $e^{t\cdot (\zeta_X-\cu)|_{\tt=t}} \in \cK_t$.

\item{\bf  The $\cA$-case} cannot be addressed by verbatim the same argument for various reasons.
  For example, the condition  $\di_t f_t=\xi_X (f_t)+\xi_Y(f_t)\in T_{\cA_t}f$ does
 not imply  $(\di_t-\xi_X-\xi_Y)^jf_t=0$ for $j\ge2$, unless $\xi_Y\in (y)\cdot T_{\cL_t}$.
   Even worse, one cannot present $e^{-\tt\cdot \di_t}e^{\tt\cdot (\di_t-\xi_X-\xi_Y)}$ as  $ e^{\tt\cdot (\zeta_X+\zeta_Y)}$
    with $\zeta_X\in Der_{\k[[t,\tt]]}(R_{X,t,\tt})$ and $\zeta_Y\in Der_{\k[t,\tt]]}(R_{Y,t,\tt})$.
 In fact,  in the BCH   formula one gets $[\di_t ,\xi_Y]\not\in T_{\cL_t}$, unless $\xi_Y=\sum l_i(y)\di_{y_i}$, where $\{l_i(y)\}$ are homogeneous
  linear forms (that depend on $t$).

Here is the adjusted argument. Starting from $\di_t f_t=\xi_X (f_t)+\xi_Y(f_t)$, rewrite it as $(\di_t-\xi_X) f_t=\xi_Y(f_t)$.
 Here $\di_t-\xi_X\in Der_\k(R_{X,t})$.
 By the algebraic Thom-Levine property,  \cite[\S3.6]{Kerner.Group.Orbits}, we get: $e^{\tt(\di_t-\xi_X)}f_t=e^{\tt\xi_Y}f_t$.
 (Note that both $\tt(\di_t-\xi_X)$ and $\tt\xi_Y$ are nilpotent derivations for the filtration $(x,t,\tt)^\bullet\sset R_{X,t,\tt}$.)
 Thus $f_t=e^{ -\tt(\di_t+\xi_X)}\circ e^{\tt\xi_Y}f_t$. As in the $\cR$-case we get:
 \beq\label{Eq.Proof.Trivial.vs.Infinit.Triv.A.case}
f_{t+\tt}= e^{\tt\cdot \di_t}\circ e^{ -\tt(\di_t-\xi_X)}\circ e^{\tt\cdot\xi_Y}f_t= e^{\tt\cdot \zeta_X}\circ e^{\tt\cdot\xi_Y}f_t,\quad\quad
 \text{for some } \quad \zeta\in Der_{\k[[t,\tt]]}(R_{X,t,\tt}).
 \eeq
 The identity
 $f_{t+\tt}=  e^{\tt\cdot \zeta_X}\circ e^{\tt\cdot\xi_Y}f_t$
is a power series in $\tt,t$, and does not involve any $t,\tt$-derivatives. Substituting $\tt=-t$ we get
 $f_o=e^{t\cdot \zeta|_{\tt=t}}\circ e^{t\cdot\xi_Y } f_t$, i.e. the trivialization.
\eei

\

We have constructed the trivializations in the formal case, $R_{X,t}=\quots{\k[[x,t]]}{J}$. In the analytic case, $R_{X,t}=\quots{\k\{x,t\}}{J}$,
  it is enough to remark that the BCH  formula gives analytic derivations, $\zeta_X$,  $\cu$, see \cite[\S2.1.ix]{Kerner.Group.Orbits}.
  (Note again, that the derivations $\tt(\di_t-\xi_X)$, $\tt\xi_Y$ are filtration-nilpotent.)
   Thus the constructed trivializations for $\cR,\cK,\cA$ are analytic.

For $\cG=\cR,\cK$ and $R_{X,t}=\quots{\k\{x,t\}}{J}$ or $R_{X,t}=\quots{\k\bl x,t\br}{J}$, one  uses the Artin approximation.
 In fact, the condition  to resolve is:  $f_o(\Phi_{X,t}(x))=f_t(x)$, resp.  $U\cdot f_o(\Phi_{X,t}(x))=f_t(x)$. Both are implicit function equations.
 Thus the formal solution implies an analytic/algebraic solution.

The $\cA$-case with $f_t$ finitely-$\cA$-determined is treated below.

\setcounter{enumi}{1}
\item {\bf The general case, $\mathbf{ char(\k)\ge0}$.} Expand the unfolding, $f_t=f_o+\sum_{j\ge d}t^j f_j$, here $f_j\in \RpX$, $d\ge 1$, and $f_d\neq 0$. We construct the $\cG_t$-trivialization
 iteratively, at each step applying the  $\cG_t$-transformation to increase $d$.
 By lemma \ref{Thm.Unfolding.Properties.Invariant.under.G.action} the initial condition $\di_t f_t\in T_{\cG_t}f_t$ is preserved under   $\cG_t$-transformations, i.e.
 $\di_t (g\cdot f_t)\in T_{\cG_t}(g\cdot f_t)$ for any $g\in \cG_t$.

  Let $f_t=f_o+t^d f_d+(t)^{d+1}$, with $f_o,f_d\in \RpX\sset R^{\oplus p}_{X,t}$, $d\ge1$.
 Then $\di_t f_t=d t^{d-1}f_d+(t)^d=\xi_t(f_t)=\xi_t(f_o)+(t)^d$, for some $\xi_t\in T_{\cG_t}$.

  We claim: $f_d= \xi(f_o)$, for some  $\xi\in T_\cG$.
  Indeed, if $char(\k)\nmid d$ then, just take the $t^{d-1}$-coefficient
 of $\xi_t$. If  $char(\k)\mid d$ then  $f_d=\xi(f_o)\in T_\cG f_o$ by our assumption of  $\cG$-separability.

 Now we take $g\in \cG_t$ corresponding to $\xi$ (using the $jet_0$-assumption for $char(\k)>0$), as follows.
\bei
\item
({\bf $\cR$-case}) We have $f_d=\xi (f_o)$.  Extend the map $x\to x-t^d\xi(x)$ to an automorphism $\Phi_{X,t}$.
 Then $\Phi_{X,t}(f_t(x))=f_t(x-t^d\xi(x)+(t)^{d+1})=f_o(x)+(t)^{d+1}$.
\item
({\bf $\cK$-case}) We have $f_d=u\cdot f_o+\xi (f_o)$.   Then $(\one-t^d\cdot u)\cdot \Phi_{X,t}(f_t)= f_o(x)+(t)^{d+1}$.
\item
({\bf $\cA$-case}) We have $f_d=\xi_Y(f_o)+\xi_X (f_o)$.   Then $(\one-t^d\cdot\xi_Y)\circ f_t\circ  \Phi_{X,t} = f_o(x)+(t)^{d+1}$.
 Note that $\one-t^d\cdot\xi_Y\in \cL_t$.
\eei

Thus for each $d\ge1$ we have $g_d\in \cG_t$ ensuring the transition $f_o+(t)^d\stackrel{g_d}{\rightsquigarrow}f_o+(t)^{d+1}$.
  This $g_d$ is of the form $\one+(t)^d$. Therefore
   the (formal) limit $\hg:=\lim(g_d\cdots g_1)\in \hat{\cG}_t$   exists, and it  satisfies: $\hg(f_t)=f_o$.

We have constructed the formal trivialization, $\hg\cdot f_t=f_o$ for $\hg\in \hat{\cG}_t$.
 For the groups $\cR$, $\cK$ one applies the Artin approximation, as before.

For the $\cA$-case, we have proved: $F\stackrel{\cA_t}{\sim}(f_o+(t)^d,t)$ for any $d\gg1$.
By the finite-$\cA_t$-determinacy of $f_t$ we get  $ f_o+(t)^d\stackrel{(\Phi_{X,t},\Phi_{Y,t})}{\sim}f_o$. Here $\Phi_{X,t}$ is $t$-linear.
 Therefore the pair $(\Phi_{X,t},\Phi_{Y,t})$ induces the final trivialization of the unfolding,  $F\stackrel{\cA_t}{\sim}(f_o,t)$.
\epr\eee

\bex
If $char(\k)=0$ then all the  unfoldings are $\cG$-separable.
\bei
\item For $\cG\in \cR,\cK$ we get:  $F$ is trivial iff it is infinitesimally  trivial.

For $R_X=\C\{x\}$ and $\cG=\cR,\cK$ this gives the classical result, see e.g. Theorem 2.22 (page 126) of \cite{Gr.Lo.Sh}.
 For $R_X=\k[[x]]$,  $char(\k)=0$ this gives Theorem 3.1.11 of \cite{Boubakri.PhD}.
\item For $\cG=\cA$ and $R_X=\quots{\k\{x\}}{J}$ we also get:  $F$ is trivial iff it is infinitesimally  trivial.

For $R_X=\C\{x\}$    this gives the classical   Thom-Levine theorem, see pg. 88 of \cite{Mond-Nuno}, for   $C^\infty(\R^n,o)$
 see \cite[pg. 176]{Martinet.1982}
\eei
Note that $f$  is not assumed   to have an isolated critical/singular/instability point.
\eex

\beR\label{Re.Trivializations.With.Constraints}
\bee[\bf i.]
\item      The trivializations  in part 2 of the theorem involve the automorphism  $\Phi_{X,t}\in Aut_\K(R_{X,t})$.
      This can be of the form $\Phi_{X,t}(x)=x+h(t)$, for $h(t)\in \K$, i.e. it does not preserve the origin of $X$.
 And similarly for the automorphism of the target, $\Phi_{Y,t}\in Aut_\K(R_{Y,t})$.
\item
      One can impose
       various restrictions on the trivializations, accordingly modifying the statements.
 For example, in the same way one gets: $\di_t f_t\in T_{\cG^{(j)}}f_t$ iff $f_t$ is $\cG^{(j)}$-trivializable.
\eee
\eeR

We give a corollary of Mather-Yau/Gaffney-Hauser type.
\bcor\label{Thm.Corol.of.M.Y.G.H.type}
Let   $\cG=\cK $ with $R_{X,t}\in  \quots{\k[[x,t]]}{J}, \quots{\k\bl x,t \br}{J},    \quots{\k\{ x,t \}}{J}$, see \S\ref{Sec.Notations.Rings.Germs}.i
  or $\cG=\cA$ with   $R_{X,t}=  \quots{\k[[x,t]]}{J}$.
 For $char(\k)>0$   suppose the unfolding $F(x,t)=(f_t(x),t)$ is  $\cG$-separable.
  If $T_{\cG_t}f_t=T_{\cG_t} f_o\sseteq \RpXt$ then $F$ is $\cG$-trivial.
\ecor
\bpr
We have: $f_t\in T_{\cG_t} f_t=T_{\cG_t} f_o$. Expand in $t$-powers, $f_t(x)=\sum_{d\ge0}t^d f_d(x)$. Then one gets:
 $f_d\in T_{\cG} f_o$ for each $d\ge0$. Therefore at the formal level (for $R_{X,t}=\quots{\k[[x,t]]}{J}$) one has:
  $\di_t f_t\in  T_{\cG_t} f_o=T_{\cG_t} f_t$. By part 2 of theorem \ref{Thm.Unfolding.Triviality.Local} we get: $F$ is formally $\cG$-trivial, i.e.
   $f_t\in \chG_t f_o$.

 For  $R_{X,t}=\quots{\k\{ x,t \}}{J},\quots{\k\bl x,t \br}{J}$ and $\cG=\cK$ one applies the Artin approximation, \S\ref{Sec.Notations.Rings.Germs}.
\epr

\subsection{An application: factorization of space-germs and map-germs}\label{Sec.Triviality.of.Unfolding.Factorization.of.Schemes.Maps}

The map of spaces $X\stackrel{f}{\to}Y$ defines the germ of (not locally trivial) fibration $\cF_f$, by $X=\coprod_{y\in Y} f^{-1}(y) $.

 Two fibrations are called equivalent, $\cF_f \sim \cF_{\tf}$ if there exist automorphisms $\Phi_X\circlearrowright X $, $\Phi_Y\circlearrowright Y$
  satisfying $\Phi_X(f^{-1}(y))=f^{-1}(\Phi_Y(y))$. Therefore one gets:  $f\stackrel{\cA}{\sim}\tf$ iff  $\cF_f \sim \cF_{\tf}$.

 Recall the classical factorization statement. Suppose a $\C$-analytic germ admits $r$ vector fields that are linearly independent at the origin.
  Then $(X,x)\cong(\C^r,o)\times(\tX,\tx)$. See e.g.  \cite[\S2.12]{Fischer}.

We establish (in zero characteristic) the $\cR$, $\cK$, $\cA$-versions of this statement.

\bed Let $(R_X,\cm)$ be a local $\k$-algebra, with $\k\cong \quots{R_X}{\cm}$.
\bee
\item The value of the derivation $\xi\in Der_X$ at the origin is the $\k$-linear map $\quots{\cm}{\cm^2}\stackrel{\xi|_o}{\to}\quots{R_X}{\cm}$.
 (Thus $\xi|_o\neq 0$ iff $\xi(\cm)\not\sseteq \cm$.)
\item The set of values of vector fields at the origin is the $\k$-vector subspace
  $Der_X|_o:=Span_\k(\xi|_o, \ \xi\in Der_X)\sset Hom_\k(\quots{\cm}{\cm^2},\quots{R}{\cm})$.
 Thus $dim_\k(Der_X|_o)$ is the maximal number of derivations linearly independent at $o\in X$.
\eee
\eed

We observe:
\bei
\item $Der_X|_o$ is invariant under $Aut^{(1)}_X$-transformations, i.e.
 is preserved by the coordinate changes with unit linear part.
\item $dim(Der_X|_o)$ is preserved by any coordinate change, $Aut_X$.
\eei

\bthe\label{Thm.Factorization.of.Spaces.Maps}
  Let $S$ be one of $\k[[x]]$, $\k\{x\}$, $\k\bl x\br$, where $char(\k)=0$. Let $R_X=\quots{S}{J}$ with $J\sseteq (x)^2$.
Suppose  $r:=dim(Der_X|_o)>0$.
\bee
\item Then $X\cong \tX\times(\k^r,o)$, i.e. the
 corresponding $\k$-algebras are isomorphic,
\[
\quots{\k[[x]]}{J}\cong \tR[[z]],\quad
\quots{\k\{x\}}{J}\cong \tR\{z\},
 \quad  \quots{\k\bl x\br }{J}\cong \tR\bl z\br.
 \]
 Here $\tR=\quots{\k[[\tx]]}{J_\tx}$, resp. $\quots{\k\{\tx\}}{J_\tx}$, resp. $\quots{\k\bl x\br }{J_\tx}$, with
  $\tx=(\tx_1,\dots,\tx_{n-r})$, $J_\tx\sseteq(\tx)^2$, and $z=(z_1,\dots,z_r)$.

\item
Moreover, take a map $X\stackrel{f}{\to}(\k^p,o)$ and the filtration $\cm^\bullet\sset R_X$.
\bei
\item ($\pmb \cR $) If $T_\cR f=T_{\cR^{(0)}}f$ then $f$ factorizes as $X\to \tX\stackrel{f|_{z=o}}{\to}(\k^p,o)$.
 Namely, $f$ is $\cR$-equivalent to a pullback of $f|_{z=o}$.
\item ($\pmb\cK$) If $T_\cK f=T_{\cK^{(0)}}f$ then the subscheme $V(f)\sset X$ factorizes into $V(f|_{z=o})\times(\k^r,o)\sset \tX\times(\k^r,o)$.
  Namely, $f$ is $\cK$-equivalent to a pullback of $f|_{z=o}:\tX\to (\k^p,o)$.
\item ($\pmb\cA$) Suppose $R_X=\quots{\k\{x\}}{J}$ or $f$ is $\cA$-finitely determined.
 If $T_\cA f=T_{\cR^{(0)}}f+T_{\cL}f$ then the fibration factorizes, $(X,\cF_f)\cong(\tX,\cF_{f|_{z=o}})\times(\k^r,o)$.
 Namely, $f$ is $\cA$-equivalent to a pull-back of $f|_{z=o}$.
\eei
\eee
\ethe
\bpr
\bee
\item Take a derivation $\xi\in Der_X$ that does not vanish at the origin, i.e. $\xi|_o\neq 0$, i.e. $\xi(\cm)\not\sseteq \cm$. Applying
 a coordinate change, i.e. $Aut_X$, one can assume $\xi(x_n)=1$. Therefore we have the
   presentation $\xi=\frac{\di}{\di x_n}-\sum^{n-1}_{i=1}a_i(x)\frac{\di}{\di x_i}$. Fix some generators $J=(q_1,\dots,q_l)\sset S$. Then for the vector $q\in S^{\oplus l}$
     one has $\di_{x_n} q\in Jac_{x_1,\dots,x_{n-1}}(q)+(q_1,\dots,q_l)\cdot S^{\oplus l}$.
      By theorem \ref{Thm.Unfolding.Triviality.Local} $q(x_1,\dots,x_n)$ is a $\cK$-trivial unfolding of $q(x_1,\dots,x_{n-1},0)$.
      Namely: $q(x_1,\dots,x_n)\stackrel{GL(l,S)\rtimes Aut_\k(S)}{\sim} q(x_1,\dots,x_{n-1},0)\in S^{\oplus l}$.
       Therefore $R_X\cong  (\quots{\k[[\tx]]}{J_\tx})[[z]]$, and similarly for the analytic/algebraic cases.

In the case $r=1$ this proves the statement. For $r>1$ one starts from the rectified germ $\tX\times (\k^1,o)$.
 Then $Der_{X\times (\k^1,o)}\ni \frac{\di}{\di x_n}$, this vector field does not vanish at $o$. By the assumption there exist
  other derivations $\xi_1, \dots ,\xi_{r-1}$, linearly independent of $\frac{\di}{\di x_n}$ at the origin. One can assume:
   $\xi_j(x_n)=0$ for $j=1\dots r$. And then repeat this procedure for $\tX$.

\item
We start from the factorized germ, $\tX\times (\k^r,o)$. Note that the condition $T_\cG f=T_{\cG^{(0)}}f$
 is preserved under the $Aut_X$-coordinate change, see lemma \ref{Thm.Unfolding.Properties.Invariant.under.G.action}.
\bee
\item[]($ \pmb \cR $) The assumption $T_\cR f=T_{\cR^{(0)}}f$ gives:
\beq
Span_{R_X}[\di_{z_1}f,\dots,\di_{z_r}f]\sseteq (\tx,z)\cdot Span_{R_X}[\di_{z_1}f,\dots,\di_{z_r}f]+Der^{(0)}_\tX(f).
\eeq
   By Nakayama we get:
$Span_{R_X}[\di_{z_1}f,\dots,\di_{z_r}f]\sseteq  Der^{(0)}_\tX(f)$. By theorem \ref{Thm.Unfolding.Triviality.Local} the unfolding $f(\tx,z)$  of $f(\tx,o)$ is $\cR$-trivial.
 Its trivialization is a $z$-linear automorphism of $R_X$. Thus it preserves the  factorized form $\tX\times (\k^r,o)$ and gives the statement.
\item[]($ \pmb \cK $) The assumption $T_\cK f=T_{\cK^{(0)}}f$ is now:
\beq
Span_{R_X}[\di_{z_1}f,\dots,\di_{z_r}f]\sseteq (z,\tx)\cdot Span_{R_X}[\di_{z_1}f,\dots,\di_{z_r}f]+Der^{(0)}_\tX(f)+(f)\cdot \RpX.
\eeq
 By Nakayama we get:
$Span_{R_X}[\di_{z_1}f,\dots,\di_{z_r}f]\sseteq  Der^{(0)}_\tX(f)+(f)\cdot \RpX$. By theorem \ref{Thm.Unfolding.Triviality.Local} the unfolding $f(\tx,z)$  of $f(\tx,o)$
 is $\cK$-trivial.
 Its trivialization preserves the  factorized form $\tX\times (\k^r,o)$ and gives the statement.
 \item[]($\pmb\cA $) The assumption is now:
 \beq
 Span_{R_X}[\di_{z_1}f,\dots,\di_{z_r}f]\sseteq (z,\tx)\cdot Span_{R_X}[\di_{z_1}f,\dots,\di_{z_r}f]+Der^{(0)}_\tX(f)+T_{\cL^{(0)}} f.
 \eeq
 By Nakayama over $R_X$, see \S\ref{Sec.Notations.Rings.Germs}.iv, we get:
$Span_k[\di_{z_1}f,\dots,\di_{z_r}f]\sseteq  Der^{(0)}_\tX(f)+T_{\cL^{(0)}} f$. By theorem \ref{Thm.Unfolding.Triviality.Local} the unfolding $f(\tx,z)$  of $f(\tx,o)$
 is $\cA$-trivial.
 Its trivialization preserves the  factorized form $\tX\times (\k^r,o)$ and gives the statement.
\epr\eee
\eee

In part 2 of this theorem we can take the filtration $(\tx,z_1,\dots,z_j)^\bullet\sset R_{\tX\times(\k^r,o)}$ for some $j\le r$. We get (for $\cR$):
 if $T_\cR f=T_{\cR^{(0)}} f$ then $f$ factorizes as $X\to \tX\times(\k^{r-j},o)\stackrel{f|_{z=o}}{\to}(\k^p,o)$.
  (And similarly in the $\cK,\cA$-cases.) The proof is the same.
We state this explicitly.
\bcor\label{Thm.Geometric.Lemma.of.Unfolding}
If $(\di_{z_1}-\xi_{\tx,z_2,\dots,z_r})f=0$, where $\xi_{\tx,z_2,\dots,z_r}\in Der^{(0)}_\tX+\cm\cdot Span_{R_X}[\di_{z_2},\dots,\di_{z_r}]$,
 then $f$ is $\cR$-equivalent to the pullback of a map  $\tX\to(\k^p,o)$.
\ecor
 This is the classical ``geometric lemma of deformations", e.g. \cite[pg.192]{Martinet.1982} (for $C^\infty$-case).

\subsection{The rank of a map}\label{Sec.Triviality.of.Unfolding.Rank.of.Map}
Given a map $F\in \cm\cdot\RpX$ we evaluate its tangent image  $T_\cR F\sseteq \RpX$ at the origin, i.e.
 take the image $\quots{R_X}{\cm}\cdot T_\cR F\sseteq (\quots{R_X}{\cm})^{\oplus p}\cong \k^{\oplus p}$.
 This is a $\k$-vector subspace.
\bed
The rank of $F$ is the dimension $dim_\k(\quots{R_X}{\cm}\cdot T_\cR F)$.
\eed
If $rank(F)=r$ then, in particular, $dim_\k(Der_X|_o)\ge r$. Thus, (assuming $char(\k)=0$) by theorem \ref{Thm.Factorization.of.Spaces.Maps} we get:
  $X\cong \tX\times(\k^r_u,o)$.
 \bel\label{Thm.prenormal.form.of.a.map}
 \bee
 \item
Take a map $F\in \cm\cdot R^{\oplus(p+r)}_X$ of rank $r$. Suppose either $char(\k)=0$ or $J=0$.
Then $F\stackrel{\cA}{\sim}(f(\tx)+h(\tx,u) ,u)$, where $u=(u_1,\dots,u_r)$ and $f(\tx)\in (\tx)^2\cdot R^{\oplus p}_\tX$ and
$h(\tx,u)\in (u)\cdot(\tx)\cdot \RpX$.
\item Here the $\cA$-type of $f$ is not uniquely determined due to the residual $\cA$-equivalence:
\[(f+h(x,u),u)\sim(f+h(x,u(\tu+q(f))),\tu)\]
 for any element $q(y)\in (y)\cdot \RpY$ and the corresponding parameter change $\tu=u+q(f)-q(f+h(x,u))$.
\eee
\eel
In particular, part 2 gives:  $(f+h(x,u),u)\sim(f+h(x,u(q(f)))+(x)(\tu),\tu)$.
\\\bpr
\bee
\item Using the subgroup $GL(p,\k)<\cL$ we can ensure: the last $r$ components of $F$ are of order $1$ (and their linear parts are independent), while
the other components of $F$ are of order$\ge2$.  By an $\cR$-transformation on the last $r$ components, we get  $F\sim ((u,x)^2,u)$.
 Then by $\cL$ we achieve the claimed form $(f(\tx)+(\tx)\cdot (u),u)$.

\item The relation $\tu=\tu(u)$ is invertible. Rewrite it as $u(\tu)=\tu-q(f)+q(f+h(x,u))$. Then:
\beq
(f+h(x,u),u)=(f+h(x,u(\tu)),u(\tu))\stackrel{\cL}{\sim}(f+h(x,u(\tu)), \tu-q(f))\stackrel{\cR}{\sim}(f+h(x,u(\tu+q(f))), \tu).
\quad\quad \epr
\eeq
\eee

\bex\bee[\bf i.]
\item Suppose $p=1$, i.e. the map is of corank=1. By Morse lemma, \cite[\S4]{Kerner.Group.Orbits},
 $F\sim(f(\tx)+Q_2(z)+h(\tx,u),u)$, where $f(\tx)\in (\tx)^3$, while $Q_2(z)$ is a non-degenerate
 quadratic form.
\item Suppose $h(x,u)=h\cdot u$, where $h\in Mat_{p\times r}(x)$. Take $q(y):=q\cdot y$, where $q\in Mat_{r\times p}(\k)$. We get
 $\tu=u-q\cdot h\cdot u$. Thus $u=[\one-q\cdot h]^{-1}\cdot \tu$. Then part 2 gives:
\[
(f+h\cdot u,u)\stackrel{ }{\sim}(f+h\cdot [\one-q\cdot h]^{-1}(\tu+q\cdot f),\tu)=
([\one-h\cdot q]^{-1}f+ h\cdot [\one-q\cdot h]^{-1}\tu,\tu).
\]
In particular, if ($)\cdot h$ spans the vector space $(x)\cdot T^1_\cK f$, then $[\one-h\cdot q]^{-1}f-f$ spans $(x)(f)\cdot \RpX\ mod T_\cK f$. Then also
 $h\cdot [\one-q\cdot h]^{-1} $ spans  $(x)\cdot T^1_\cK f$. Thus any two $\cA$-stable maps that are $\cK$-equivalent are also $\cA$-equivalent.
\eee
\eex

\subsection{Trivializing unfoldings over non-local bases}\label{Sec.Triviality.of.Unfolding.Nonlocal.Base}
Suppose an unfolding is  infinitesimally trivial.
Theorem \ref{Thm.Unfolding.Triviality.Local} ensures the local triviality of the unfolding over the base $ (\k^r_t,o)$.
In the analytic case, $R_X=\quots{\k\{x\}}{J}$, one gets the trivialization over a small ball, for $\|t\|\ll1$.
 In the henselian case, $R_X=\quots{\k\bl x\br}{J}$, one gets the trivialization over an \'etale neighborhood.

Take an unfolding with the affine base, $F=(f_t(x),t)$ with $t\in \k^r$. The ring $R_{X,t}$ is now
 one of $\quots{\k[t][[x]]}{J}$,  $\quots{\k[t]\{ x\}}{J}$,  $\quots{\k[t]\bl x\br}{J}$, i.e. power series in $x$ whose coefficients are polynomials in $t$.
Take the corresponding group $\cG_t\in \cR_t,\cK_t,\cA_t$. These transformations are $t$-global.

Accordingly we want the $t$-global trivialization criteria.  These are obstructed in two ways.
 \bei
 \item  The $\cG_t$-transformations do not preserve the origin. Therefore, e.g. for $R_X=\C\{x\}$, when taking  non-small $t$, one runs out of the
  ball of convergence in $\C^n$.
 Therefore one should restricts to the subgroup $\cG^{(o)}_t\le\cG_t$ of elements that
   preserve the origins of $X$ and $(\k^p,o)$.
\item
 The filtered-nilpotent vector fields can be ``integrated", i.e. there exists a map $T_{\cG^{(1)}}\to \cG^{(1)}$ that approximates the classical exponential,
  \cite[\S3]{BGK.20}. But not every element $\xi\in T_{\cG^{(o)}}$ induces an element   $g_\xi\in \cG^{(o)}$ satisfying: $ord(g_\xi-Id-\xi)>ord(\xi)$.
 For example, for $tx\di_x \in T_{\cR^{(o)}}$ the corresponding $\cR$-transformation would be of the type $x\to x+tx+(x)^2$. But this is not invertible for $t=1$.
\eei

Therefore below we restrict to the subgroup $\cG^{(1)}_t<\cG_t$.
 For the subgroup $\cG^{(0)}_t<\cG_t$ our results are weaker, we can only compare particular fibres of the family $\{f_t\}$.

\subsubsection{Trivialization for the subgroup $\cG^{(1)}< \cG$}
\bel
Let $R_{X,t}=\quots{\k[t][[x]]}{J}$,   where $\k$ is a field of zero characteristic.
 Let $\cG\in \cR,\cK,\cA$.  Any infinitesimally $\cG^{(1)}$-trivial unfolding is $t$-globally $\cG^{(1)}$-trivial.
\eel
\bpr
We use the trivializations constructed in the proof of theorem \ref{Thm.Unfolding.Triviality.Local}, equations
  \eqref{Eq.Proof.Trivial.vs.Infinit.Triv.R.case},  \eqref{Eq.Proof.Trivial.vs.Infinit.Triv.K.case},  \eqref{Eq.Proof.Trivial.vs.Infinit.Triv.A.case}.
   We should only verify: the derivation $\zeta_X$ constructed via
   $e^{-\tt\di_t}\circ e^{\tt(\di_t+\xi_X)}=e^{\tt\zeta_X   } $ belongs to $T_{\cR^{(1)}_t}$.
 In particular, in its $(t,x)$-expansion all the coefficients of $x$-monomials should be polynomials in $t$ (rather than just power series).

    Indeed, by the BCH formula, \cite[\S2.1]{Kerner.Group.Orbits}, $\zeta_X=\sum \tt^l p_l(\di_t,\di_t -\xi_X)$, where $p_l$ is a homogeneous polynomial expressible via repeated
     commutators of $\di_t$ and $\xi_X$. Observe:
     \bei
     \item $[T_{\cR^{(i)}_t},T_{\cR^{(j)}_t}] \sseteq  T_{\cR^{(i+j)}_t}$.
\item For each $\eta\!\in\! T_{\cR^{(j)}_t}$ we have $[\di_t,[\di_t,\dots,[\di_t,\eta]]]\!\in\! T_{\cR^{(j+1)}_t}$, for the $\di_t$-commutator
 repeated sufficiently many times.
 This holds because the coefficient of every $x$-monomial is a polynomial in $t$.
     \eei
Therefore every summand that appears in $\sum \tt^l p_l(\di_t,\di_t -\xi_X)$ belongs to $T_{\cR^{(1)}_t}$, and moreover,
  $p_l(\di_t,\di_t -\xi_X)\in T_{\cR^{(j)}_t}$ for each $j$ and a corresponding $l\gg1$. Altogether, $\zeta_X\in T_{\cR^{(1)}_t}$.
\epr

\subsubsection{Comparison of fibres for the subgroup $\cG^{(0)}\le \cG$}

\bprop
Let $R_X$ be one of $\quots{\k[[x]]}{J}$, $\quots{\k\{x\}}{J}$, $\quots{\k\bl x\br}{J}$, with $\k$ either $\R$  or alg.closed field of zero characteristic.
For $\cG=\cA$ take $R_X=\quots{\k[[x]]}{J}$.  Take an unfolding  $F=(f_t(x),t)\in R_X^{p+r}[t]$. Suppose $\sqrt{I}=\cm$.
 If $ \di_t f_t\in T_{\cG^{(0)}_t}f_t$ then $f_o\stackrel{\cG^{(0)}}{\sim}f_{t_o}$
  for every $t_o\in \k^r$.
\eprop
\bpr
\bee[\bf Step 1.]
\item We prove the statement for $\k\in \R,\C$, and the ring $\quots{\k\{x\}}{J}$. Take a path   $o\rightsquigarrow t_o$ (in $\k$).
 Along this path (at every point on it)
  we have  $\di_t f_t\in T_{\cG^{(0)}_t}f_t$. By theorem \ref{Thm.Unfolding.Triviality.Local} we get: $F$ is
  $\cG$-trivializable on a small ball $Ball(t')$. More precisely, the trivialization is done by a representative of an element of   $\cG$.
 In fact $F$ is $\cG^{(0)}$-trivializable, see remark \ref{Re.Trivializations.With.Constraints}.ii.

  Cover the path  $o\rightsquigarrow t_o$ by such balls, and choose a finite subcover. Then we can pass from $f_0$ to $f_{t_o}$ by a finite number of $\cG^{0}$-elements,
   or rather by  their representatives. All these representatives preserve the origin, and are defined on the same (small) balls in $\k^n,\k^p$.
    Thus their product is defined, and we get: $f_{o}\stackrel{\cG^{(0)}}{\sim}f_{t_o}$.

\item (The general case, $\k=\bar\k$.)
 Let $\di_t f_t=\xi(f_t)$ for $\xi\in T_{\cG^{(0)}_t}$ and $R_{X,t}=\quots{\k[t][[x]]}{J}$. Pass to the finite jets $\quots{\k[t][[x]]}{J+(x)^d}$.
 Note that $\cR^{(o)}$ and $T_{\cR^{(0)}_t}$ act on this quotient. We still have the condition
  $\di_t [f_t]=[\xi]([f_t])$.
  Let $\{C_\bullet\}$ be the (finite) collection of $\k$-coefficients in the Taylor   expansion of $[f_t]$ and $[\xi]$.
   Then $\Q\{C_\bullet\}\supseteq\Q$ is   a finite extension of fields. Therefore one can embed $j:\Q\{C_\bullet\} \hookrightarrow \C$.

We still have   $\di_t j(f_t)=j(\xi)\cdot j(f_t)$, now over the ring $\quots{\C[t][[x]]}{j(J)+(x)^d}\cong \quots{\C[t]\{x\}}{j(J)+(x)^d}$.
 Therefore, by part 1, we get  $j[f_{o}]\stackrel{\cG^{(0)}_\C}{\sim}j[f_{t_o}]$.
 Consider this as the equivalence $[f_{o}]\stackrel{\cG^{(0)}}{\sim}[f_{t_o}]$ over some field extension of $\k$.
 By lemma \ref{Thm.Equivalence.Change.base.field} we get the equivalence over $\k=\bar\k$. (Now the polynomial system to resolve is finite,
  thus we do not need the assumption ``$\bar\k$ is uncountable".)

Returning to $\quots{\k[t][[x]]}{J}$, we get (for every $d\ge 1$) an element $g_d\in \cG^{(0)}$ satisfying  $f_{o}\stackrel{\cG^{(0)}}{\sim}f_{t_o}\ mod(x)^d$.
 Applying this iteratively, we can present this via the product,
 $f_{t_o}=g_d\cdots g_1\cdot f_o$, with $g_j\in \cG^{(j)}$. This product converges formally. Therefore $f_{o}{\sim}f_{t_o}$ are $\cG^{(0)}$-equivalent.

For $\cR,\cK$ one applies the Artin approximation.
\epr
\eee

\beR
\bee[\bf i.]
\item In the first part of this proof we needed the normed field $\k$ to be connected and locally compact. This forces $\k$ to be $\R$ or $\C$,
 \cite{Wieslaw}.
\item The statement does not hold over fields that are not $\R$-closed. Even the $t$-local version (over $\k[t]_{(t)}[[x]]$) fails.
 For example, $f_t:=(1+t)x^d\in \k[t]_{(t)}[[x]]$ satisfies $\di_t f_t\in T_{\cR^{(0)}_t} f_t$.
 But to achieve $f_o\sim f_{t}$, at least for $|t|\ll1$, one needs the property  $\sqrt[d]{1+t}\in \k$ for $|t|\ll1$.
\eee
\eeR

\section{The pre-normal form of unfolding and versality}\label{Sec.Versality}
\subsection{}\label{Sec.Versality.Prenormal.Form}
 Let $R_X=\quots{\k[[x]]}{J}$, see \S\ref{Sec.Notations.Rings.Germs}, and fix a map $f_o\in \cm\cdot \RpX$, not necessarily $\cG$-finite.
 We get the $\k$-vector space $T^1_\cG f_o:=\quots{\RpX}{T_{\cG}f_o}$, possibly of infinite dimension.
  Fix (any) elements $\{v_\bullet\}\sset \RpX$ whose images formally generate the $\k$-vector space $T^1_\cG f_o$.
 Namely, any element of $T^1_\cG f_o$ is
 presentable as $\sum c_\bullet v_\bullet$, with $c_\bullet\in \k$, and this sum converges in $R_X$.
 For $\cG=\cA$ we always assume  $\{v_\bullet\}\sset \cm\cdot\RpX$.
 \bed
  The $\cG$-pre-normal form of an unfolding of $f_o$ is the unfolding $(f_o(x)+\sum a_j(t )\cdot v_j(x),t)$, where
   $a_j(t)\in (t)\cdot \k[[t]]$.
 \eed

\bel\label{Thm.Unfolding.Normal.form}
Any unfolding  $F(t,x)=(f_t(x),t)\in R^{\oplus (p+r)}_{X,t}  $ is (formally) $\cG$-equivalent to its pre-normal form.
Moreover, this $\cG$-equivalence preserves the image of $Span_\k[\di_t f_t|_{t=o}]$ inside $T^1_\cG f_o$, i.e.
\beq\label{Eq.Prenormal.Form.Unfolding}
T_{\cG}f_o+Span_\k[\di_t f_t|_{t=o}]=T_{\cG}f_o+Span_\k[\sum \di_t a_j(t)|_{t=o}\cdot v_j(x)]\sseteq\RpX.
\eeq
\eel
\bpr The transition to the pre-normal form is inductive.
Using the generators $\{v_\bullet\}$ we can present $F(t,x)=(f_o+\xi^{(d)} f_o+\sum a^{(d)}_j(t)v_j,t)$.
 Here $1\le d<\infty$ and $0\neq \xi^{(d)} f_o\in (t)^d\cdot T_{\cG_t} f_o$,   and $ a^{(d)}_\bullet (t)\in (t)\sset \k[[t]] $.
 Take a group element $g_d\in \cG_t$ of the form $g=Id- \xi^{(d)}+(t)^{d+1}$.
   (This is ensured by the $jet_0$-assumption on $R_X$, see \S\ref{Sec.Notations.Rings.Germs}.)
  Then $g_d(F)=(f_o+\xi^{(d+1)} f+\sum a^{(d+1)}_j(t)v_j,t)$, for some new power series $\{a^{(d+1)}_\bullet(t)\}$.
  Iterate this process. Note: $a^{(d+1)}_j(t)-a^{(d)}_j(t)\in (t)^d$.

  The infinite product $\lim_{d\to\infty}(g_d\cdots g_1)$ converges (formally) to an element   $g\in\cG$.
    Indeed, at $d$'th step we
   do not change the $t^i$-elements with $i<d$.

 Altogether, the unfolding $gF$ is in the pre-normal form.
Finally we verify equation \eqref{Eq.Prenormal.Form.Unfolding}.
 \bei
\item  $g_1$ changes $Span_\k[\di_t f_t|_{t=o}]$ only by an element of  $T_{\cG} f_o $. Hence $g_1$
preserves $T_{\cG} f_o +Span_\k[\di_t f_t|_{t=o}]$.
 \item For $d\ge2$ the element $g_d$ does not affect $\di_t f_t|_{t=o}$.
\epr
\eei

\beR\label{Rem.Pre-normal.form}
\bee[\bf i.]
\item The pre-normal form  is usually far from being unique. E.g. suppose $g f_o=f_o$ for some $Id\neq g\in \cG$. Then the ($\cG$-equivalent)
 unfoldings $f_t,gf_t$ have different pre-normal forms.
\item Suppose the ring is $R_X=\quots{\k\bl x\br}{J}$, the group is $\cG\in \cR,\cK$ and $f_t$ is a $\cG_t$-finite map.
  Then the pre-normal form can be achieved
 over $R_X$, not just formally. Indeed, we should resolve the condition $g\cdot f_t=f_o+\sum a_j(t)v_j$, with the unknowns $g,\{a_\bullet\}$.
  Here $g$ can depend on $x,t$,
  but $a_\bullet$ should depend on $t$ only.   This is a (finite) system of $\k\bl x,t\br$-equations.
  The lemma ensures a formal solution  $\hg,\{\ha_\bullet(t)\}$, over $\hR_{X,t}$.
Apply the nested Artin approximation   to achieve an ordinary solution, $g,\{a_\bullet(t)\}$.

\item In the $\cA$-case another presentation of pre-normal form is often useful. Suppose $f$ is $\cK$-finite.
 Fix some elements $\{v_\bullet\}\sset \cm\cdot \RpX$ that go to a basis of $\cm\cdot T^1_\cK f_o$.
 Then $\RpX=Span_{R_Y}(v_\bullet)+T_\cA f_o$, see lemma \ref{Thm.TA.vs.TK}.
  And then,    lemma \ref{Thm.Unfolding.Normal.form} gives: any unfolding of $f_o$ is formally $\cA$-equivalent to
   $(f_o +\sum a_j(t,f_o)\cdot v_j ,t)$, for some $a_j(t,y)\in (t)\sset\k[[t,y]]$.
\eee
\eeR
\bex
\bee[\bf i.]
\item Take a map  $f:(\k^n,o)\to (\k^n,o)$ of corank one. Present it as an unfolding, $f=(u,f_o(x)+(x)(u))$, here $u=(u_1,\dots,u_{n-1})$.
  Then $f\stackrel{\cA}{\sim}(u,x^d+\sum^{d-2}_{i=1}x^i a_i(u))$,  with $a_j(u)\in (u)$.
  \item
Take a map  $f:(\k^n,o)\to (\k^{n+1},o)$ of corank one. Present it as an unfolding, $f=(u,f_1(x)+(x)(u),f_2(x)+(x)(u))$.
 By $\cA$-equivalence we can assume $f_1(x)=x^{d}$, $ord(f_2)>d$, and $d \nmid ord(f_2)$.
  Then $f\stackrel{\cA}{\sim}(u,x^{d}+\sum_{j=1}^{d-2}a_j(u)x^j, f_2(x)+(x)(u))$, with $a_j(u)\in (u)$.
\eee
\eex

\subsection{Versality vs infinitesimal versality}\label{Sec.Versality.vs.Infinitesimal.Versality}
 Let $\cG=\cA$ with  $R_X=\quots{\k[[x]]}{J}$ or  $\cG\in \cR,\cK$ with $R_X$   one of $\quots{\k[[x]]}{J}$, $\quots{\k\bl x\br}{J}$,
  or  $char(\k)=0$ with $\cG\in \cR,\cK,\cA$ and  $R_X=\quots{\k\{ x\}}{J}$.
\bthe\label{Thm.Unfolding.Versality.vs.Infinites.versality}
\bee
\item Let a (finite) tuple $\{v_\bullet\}$ in $\RpX$ generate  $T^1_\cG f$.
 For $\cG\in \cR,\cK$ the unfolding  $(f_o+\sum t_j v_j ,t)$ is   $\cG$-versal.
 For $\cG=\cA$  the unfolding  $(f_o+\sum t_j v_j ,t)$ is formally $\cA$-versal.
\item
 ($\cG\in \cR,\cK$) An unfolding $F $ is   $\cG$-versal iff it is infinitesimally $\cG$-versal.
\\
($\cG=\cA$) An unfolding $F $ is   formally-$\cA$-versal iff it is infinitesimally $\cA$-versal.
\item $dim_\k T^1_\cG f$ is the minimal number of parameters in a $\cG$-versal unfolding.
\eee
\ethe
\bpr
\bee
\item Start from a(ny) unfolding $F=(f_t(x),t)\in R^{p+r}_{X,t}$. We should resolve the condition $g(f_t)=f_o+\sum a_jv_j$, with the unknowns $g\in \cG_t$
 and $a_j\in \k[[t]]$, resp. $\k\bl t\br$. The pre-normal form, lemma \ref{Thm.Unfolding.Normal.form}, ensures the formal solution. In the $\cR,\cK$-cases
 and $R_X= \quots{\k\bl x\br}{J}$ one uses remark \ref{Rem.Pre-normal.form}.

\

For  $char(\k)=0$ we give another proof for  $R_X=\quots{\k\{ x\}}{J}$, the direct generalization of \cite[Chapter XIV]{Martinet.1982}.
  Extend the ring by a new variable $\tt$ and extend the unfolding, $\tF:=(f_{t\tt},t,\tt)$, with $f_{t\tt}(x):=f_t(x)+\sum \tt_jv_j$.
 We prove: $\tF$ is induced from the unfolding $(f_o(x)+\sum t_jv_j,t)$. By lemma \ref{Thm.Unfoldings.Algebraic.Lemma} we get:
  $R^{\oplus p}_{X,t,\tt}=T_{\cG_{t\tt}}f_{t\tt}+Span_\K(\{v_\bullet\})$, where $\K=\k[[t,\tt]]$, resp. $\k\{t,\tt\}$, resp. $\k\bl t,\tt\br$.
   Therefore we have
     $\di_{\tt} f_{t\tt}\in T_{\cG_{t\tt}}f_{t\tt}+Span_\K(\{v_\bullet\})=T_{\cG_{t\tt}}f_{t\tt}+Span_\K(\di_{\tt_1}f_{t\tt},\dots,\di_{\tt_r}f_{t\tt})$.
      By corollary \ref{Thm.Geometric.Lemma.of.Unfolding} one has $f_{t\tt}\stackrel{\cR}{\sim}f_t$

Take the filtration $(\tt)^\bullet\cdot R^{\oplus p}_{X,t,\tt}$. Then $T_{\cG_{t\tt}}f_{t\tt}=T_{\cG^{(0)}_{t\tt}}f_{t\tt}$.
 Thus (theorem \ref{Thm.Factorization.of.Spaces.Maps} and corollary \ref{Thm.Geometric.Lemma.of.Unfolding}) $f_{t\tt}$ is $\cG$-equivalent to the pullback of $f_{t\tt}|_{\tt=0}=f_t$.

\item
{\bf The part $\Lleftarrow$.} Let $\K=\k[[t]]$ for $\cG=\cA$ and  $\K=\k[[t]],\k\bl t\br$ for $\cG=\cR,\cK$.

By part 1, $F$ is $\cG$-equivalent to $\tF=(f_o+\sum a_j(t)v_j,t)$, with some coefficients $a_\bullet(t)\in \K$.
 Moreover, $\tF$ is still infinitesimally-versal, by lemma \ref{Thm.Unfolding.Properties.Invariant.under.G.action}.
  Thus the elements $\di_{t_i}(\sum a_j(t)v_j)|_{t=o}$,
   $i=1,\dots,r$,  generate the vector space $T^1_\cG f_o$.  But then the map $t\to \{a_\bullet(t)\}$ has a section,
   $\{a_\bullet(t)\} \to t$. Therefore the unfolding $(f_o+\sum t_j v_j,t)$ is a pullback of $\tF$. Hence $\tF$ is versal.

\noindent{\bf The part $\Rrightarrow$.} If $F=(f_t(x),t)$ itself is versal then both $F$ is induced from  $(f_o+\sum \tt_j v_j ,\tt)$ and vice versa. And thus
   $\di_{t_1}f_t|_{t=o},\dots, \di_{t_r}f_t|_{t=o}$ generate $T^1_\cG f_o$.

   \item For any unfolding $F=(f_t(x),t)$ the number of parameters is at least $dim_\k Span_\k[\di_{t_1}f_t ,\dots, \di_{t_r}f_t]|_{t=o}$.
   And (by part 2) for a versal unfolding $dim_\k Span_\k[\di_{t_1}f_t ,\dots, \di_{t_r}f_t]|_{t=o}\ge dim T^1_\cG f_o$.
\epr
\eee

\bex
\bee[\bf i.]
\item This theorem holds also for the ring $\quots{\k\{x\}}{J}$ in positive characteristic, provided  $f_o$ is $\cG$-equivalent to a polynomial map.
 In this case one chooses $\{v_\bullet\}$ as polynomials in $x$ and works in the ring  $\quots{\k\bl x\br}{J}$.
\item
For $R_X=\R\{x\},\C\{x\},C^\infty(\R^n,o)$ this is a classical theorem.   See e.g.  \cite[Chapter XIV]{Martinet.1982} (for  $\cA$ equivalence in $C^\infty$ case),
 \cite[Theorem 9.3]{Damon84} (for $\cA,\cK$ in all cases),  pg. 143 of \cite{Mond-Nuno},
 or  Theorem 1.16, pg.238 of \cite{Gr.Lo.Sh} (for $\cK$-equivalence of $\C\{x\}$).
 We remark: the formal statement is simple, but the $\k\{x\}$, $C^\infty$-cases are non-trivial.
\eee\eex

\subsection{Fibration of $\cK$-trivial unfoldings   into $\cA$-unfoldings}\label{Sec.Versality.Fibration.of.K.unfoldings.into.A}
 Let $f_t$ be a $\cK$-trivial family. To understand the $\cA$-types appearing in this family we should take a slice transversal  to the subspace
  $T_\cA f_o\cap T_\cK f_o\sseteq T_\cK f_o$. Take the quotient $\quots{T_\cK f_o}{T_\cA f_o\cap T_\cK f_o}$. Both parts contain $T_\cR f_o$, hence the
   standard isomorphism of $R_Y$-modules gives:
   \beq\label{Eq.transversal.to.TA.inside.TK}
   \quot{T_\cK f_o}{T_\cA f_o\cap T_\cK f_o}\cong  \quot{(x)\cdot(f_o)\cdot \RpX}{T_\cA f_o\cap (x)\cdot(f_o)\cdot \RpX}\cong
     \quot{(x)\cdot(f_o)\cdot \RpX+(y)\cdot T_\cL f_o}{T_\cA f_o}.
   \eeq
 Assume  $f_o$ is $\cK$-finite, and    fix some elements $\{v_\bullet\}\sset (x)\cdot(f_o)\cdot \RpX$ whose
  images generate the $R_Y$-module $\quots{(x)\cdot(f_o)\cdot \RpX}{T_\cA f\cap (x)\cdot(f_o)\cdot \RpX}$.

Evidently any unfolding of type $(f_o(x)+(t)\cdot Span_{R_Y}\{v_\bullet\},t)$ is $\cK$-trivial.
 \bed
The $\cA$-pre-normal form of a $\cK$-trivial unfolding of the map $f_o\in (x) \cdot\RpX$ is the unfolding of type $(f_o +(t)\cdot Span_{R_Y}\{v_\bullet\},t)$.
 \eed

\bel\label{Thm.Unfolding.K-trivial.Normal.form}
Let $R_X$ be one of $\quots{\k[[x]]}{J}$, $\quots{\k\{x\}}{J}$, $\quots{\k\bl x\br}{J}$. If $char(\k)>0$ then assume $J=0$.
 Any $\cK$-trivial unfolding  $(f_t(x),t)$ of $f_o\in (x)\cdot\RpX$  is formally-$\cA_t$-equivalent to its $\cA$-pre-normal form.

Moreover, if the map $f_t \in (x)\cdot\RpXt$ is $\cA_t$-finitely determined then the unfolding $(f_t(x),t)$  is $\cA_t$-equivalent to its pre-normal form.
\eel
In this sense the unfolding $(f_o+\sum t_j v_j,t)$ is versal among all the $\cK$-trivial unfoldings.
\\\bpr
As the unfolding $(f_t(x),t)$ is $\cK$-trivial, we can assume (by an $\cR_t$-transformation): $(f_t(x))=(f_o(x))\sset R_{X,t}$.
  Then $\RpXt\ni f_t(x)\stackrel{GL(p,R_{X,t})}{\sim}f_o(x)$. Thus  after a $GL(p,\k[[t]])$-transformation, we can assume:
  \[
  f_t-f_o\in (t)\cdot (x)\cdot (f_o)\cdot \RpXt=(t)\cdot(Span_{R_{Y,t}}\{v_\bullet\}+T_\cA f_o\cap (x)\cdot(f_o)\cdot \RpX).
\]
Therefore   $  f_t-f_o-t(\xi_X(f_o)+\xi_Y(y)|_{f_o})\in (t)\cdot Span_{R_{Y,t}}\{v_\bullet\} $, for some derivations $\xi_X\in T_\cR,\xi_Y\in T_\cL$.

 If $R_X$ is regular (i.e. $J=0$) then define the coordinate change $\Phi_{X,t}\in \cR_t$ by $x\to x+t\xi_X(x)$.
 In the general case define $\Phi_{X,t}$  as the extension  $x\to x+t\xi_X(x)+(t)^2$ of the map $x\to x+t\xi_X(x)$.
 (Our assumptions ensure the $jet_0$-condition.) We remark that $\Phi_{X,t}$ does not necessarily preserve the origin of
  the source, i.e. $\Phi_{X,t}(x)\neq(x)\sset R_{X,t}$.

 Define $\Phi_{Y,t}\in \cL_t$ by $y\to y+t\xi_Y(y)$. It does not necessarily preserve the origin of the target, $\Phi_{Y,t}(y)\neq(y)\sset R_{Y,t}$.
 Altogether we get: $\Phi^{-1}_{X,t}\circ\Phi^{-1}_{Y,t}(f_t)-f_o\in (t)\cdot Span_{R_{Y,t}}\{v_\bullet\} +(t)^2\cdot T_\cA f_o.$

 Iterate this process to get the formal $\cA$-equivalence, $(f_t ,t)\sim (f_o+(t)\cdot Span_{R_{Y,t}}\{v_\bullet\},t)$.

  Suppose $f_t$ is $\cA$-finitely determined.
We have proved: $f_o+Span_{R_Y}(\{v_\bullet\})\stackrel{\cA_t}{\sim}f_t+(t)^d$ for $d\gg1$. Now invoke the finite determinacy.
\epr

\bex ($p=1$) Take a map $f:(\k^n,o)\to (\k^1,o)$, $f(x)\in (x)^2$.
 Suppose $f\in (x)^2$ is ``almost weighted homogeneous" in the following sense: $(x)\cdot f\sseteq   Jac(f)+Span_\k(f,f^2,\dots)$.
 Then in equation \eqref{Eq.transversal.to.TA.inside.TK} one get $\quots{T_\cK f_o}{T_\cA f_o}=0$.
  Therefore any $\cK$-trivial unfolding of $f$ is also (formally) $\cA$-trivial.
\eex

\subsection{Fibration of $\cK$-orbits into $\cA$-orbits}\label{Sec.Versality.fibration.of.K.orbits.into.A}
 Assume $char(\k)=0$ or $J=0$.
  Take a map of rank $r$,
 present it as $F=(f(x,u),u)\in R^{\oplus(p+r)}_{X,u}$,
 where $f(x,u)\in (x)\cdot(x,u)\in R^{\oplus p}_{X,u}$
 (see lemma \ref{Thm.prenormal.form.of.a.map}).
  Here $R_{X,u}$ is one of $\quots{\k[[x,u]]}{J}$, $\quots{\k\{x,u\}}{J}$, $\quots{\k\bl x,u\br}{J}$. Accordingly
  $R_{Y,u}$ is one of $\k[[y,u]]$, $\k\{y,u\}$, $\k\bl y,u\br$.

 Suppose the map $f_o(x)=f(x,o)\in (x)^2\cdot \RpX$ is $\cK$-finite.
 Fix some elements $\{v_\bullet\}\sset (x)\cdot\RpX$ that go
  to a basis of $(x)\cdot T^1_\cK f_o$. By the direct check: these elements go also to the basis of $(x)\cdot T^1_\cK F$.

By lemma \ref{Thm.TA.vs.TK}: $T_\cK F=T_\cR F+  T_{\cL^{(0)}} F+Span_{R_Y}\{(y)\cdot v_\bullet\}$.
 In words: the $T_\cK$ orbit decomposes into the $T_\cA$-orbit and a finitely generated complement.
 This is an ``infinitesimal" statement, we prove its local version.

  The orbit $\cK F$ splits into  $\cA$-orbits.  We consider only the orbits of $\cA$-finitely determined maps.
 \bprop\label{Thm.Orbits.K.fibers.into.A}
   Suppose $\k$ is an infinite field. Then:
\[\cK F\cap (\cA\text{-finitely determined})\sseteq \cA\Big(\{f\}+Span_{R_{Y,u}}[(y,u)\cdot \{v_\bullet\}],u\Big).
\]
 \eprop
\bpr Suppose a map $\tF\in \cK F$ is $\cA$-finitely determined. Then to bring it to the claimed form it is enough to work
 over $R_{X,u}=\quots{\k[[x,u]]}{J}$ and $R_{Y,u}=\k[[y,u]]$.
\bee[\bf Step 1.]
\item  By an $\cR$-transformation we can assume $(\tF)=(F)\sset R_{X,u}$.
 Then, by lemma \ref{Thm.prenormal.form.of.a.map}, we can assume
  $\tF=(\tf(x,u),u)\in R^{\oplus(p+r)}_{X,u}$, with $f,\tf\in (x)\cdot (x,u)\cdot \RpXu$. Here $(\tf,u)=(f,u)\sset R_{X,u}$.
 Apply $GL(p,\k[[u]])$ to get: $\tf-  f \in(x)\cdot (f,u)\cdot \RpXu$.

By lemma \ref{Thm.TA.vs.TK} we have $(x)\cdot\RpX= T_\cR f_o+  T_{\cL^{(0)}} f_o+Span_{R_Y}(\{v_\bullet\})$.
 Therefore
 \beq\label{Eq.Fibration.K.orbits.into.A.orbits}
 (x,u)\cdot R^{\oplus p}_{X,u}= T_\cR f +(y,u)\cdot T_\cL f+Span_{R_{Y,u}}[\{v_\bullet\}].
 \eeq
  Hence
 $\tf-  f \in Span_{R_{Y,u}}[(y,u)\{v_\bullet\}]+(f,u)T_\cR f +(y,u)^2\cdot T_\cL f$.

As $\tF$ is  $\cA$-finitely determined,  it is enough to
 prove: $\tF+h_d\in \cA\big(f +Span_{R_{Y,u}}[(y,u)\cdot \{v_\bullet\}],u\big)$ for some $h_N\in (x,u)^N\cdot R^{\oplus(p+r)}_{X,u}$, with $N\gg1$.
  Therefore it is enough to prove:
\beq\label{Eq.inside.proof.K.foliates.into.A.form.to.achieve}
  \tF \in \cA\big(f +Span_{R_{Y,u}}[(y,u)\cdot \{v_\bullet\}],u\big)+(x,u)^N\cdot R^{\oplus (p+r)}_{X,u} \quad \text{ for } N\gg1.
\eeq

  \item
  Start from $\tf=f+tg$, where $g\in (x)\cdot (f,u)\cdot  \RpXu$, and $t$ is an indeterminate.
The transition to the form of \eqref{Eq.inside.proof.K.foliates.into.A.form.to.achieve} is done inductively.
 By equation \eqref{Eq.Fibration.K.orbits.into.A.orbits}
  we can present $g\in Span_{R_{Y,u}}[(y,u)\cdot \{v_\bullet\}]+(f,u)^d\cdot T_\cR f+ (y,u)^{d+1}\cdot T_\cL f$ for some $d\ge1$.
   \bei
   \item
  {\bf Inductive step for the case of a regular ring, $\pmb{ R_{X,u}=\k[[x,u]]}$.}  Define the coordinate change $\Phi_X\in \cR$ by $x\to x+t(f,u)^d\cdot c\cdot\xi(x)$ and $u\to u$.
 (Here $c\in R_{X,u}$ is an unknown.)
 Note: this coordinate change is filtration-unipotent, i.e. $\Phi_X-Id$ is filtration-nilpotent for $(x,u)^\bullet\sset R_{X,u}$.

    Accordingly
    \beq
\Phi_X(f+tg)-f\in t(f,u)^d(c+c^{(g)})T_\cR f+t^2\cdot(f,u)^{d+1}\cdot(c)^2\RpXu+t^2(f,u)^{d}\cdot(c)\cdot(x)\RpXu.
    \eeq
Expand $c=\sum c_{\um,\un}x^\um\cdot u^\un$, with multi-indices $\um,\un$. Here $\{c_{\um,\un}\}\in \k[[t]]$ are unknowns. Then we get:
    \beq
\Phi_X(f+tg)-f\in t(f,u)^d\cdot \sum\nolimits_{\um,\un}\big(c_{\um,\un}+c_{\um,\un}^{(g)}+t\cdot H_{\um,\un}(\{c_*\})\big)\cdot x^\um\cdot u^\un\cdot T_\cR f+
    \eeq
\[\hspace{7cm}
   t^2(\{c_*\})^2\cdot Span_{R_Y}[(y)\cdot \{q_\bullet\}]+ t(\{c_*\})\cdot Span_{R_Y}[(y)^{d+1}\{q_\bullet\}].
   \]
In the infinite summation $\sum_{\um,\un}$ we need only the finite part  (by the finite determinacy). Thus we get a finite system of {\em polynomial} equations,
\beq\label{Eq.inside.proof.K.fibrations.system.of.IFT}
\{c_{\um,\un}+c_{\um,\un}^{(g)}+t\cdot H_{\um,\un}(\{c_*\})=0\}_{\um\un}.
\eeq
Apply $IFT_\one$ to get the solution $c_{\um,\un}(t)\in \k[[t]]$. For these coefficients the transformation $\Phi_X$ satisfies:
\beq
\Phi_X(f+tg)-f\in (x,u)^N\cdot \RpXu+Span_{R_{Y,u}}[(y)\cdot \{v_\bullet\}]+t\cdot  (y,u)^{d+1} T_\cL f.
\eeq

 Apply a (filtration-unipotent) $\cL$-transformation to eliminate the term $t\cdot  (y,u)^{d+1} T_\cL f$. We get:
\beq\label{Eq.inside.proof.K.foliates.end.of.induction}
\Phi_Y\Phi_X(f+tg)-f\in (x,u)^N\cdot \RpXu+Span_{R_{Y,u}}[(y,u)\cdot \{v_\bullet\}] +t\cdot  (y,u)^{d+1}\cdot (x)\RpXu.
\eeq
Now expand $(x)\RpXu$ as in Step 1, and iterate. In at most $N$ steps one gets equation \eqref{Eq.inside.proof.K.foliates.into.A.form.to.achieve}.

\item{\bf Inductive step for the the general case, $R_{X,u}=\quots{\k[[x,u]]}{J}$.} The previously found transformation,
  $\tilde\Phi_X: x\to x+t(f,u)^d\cdot c\cdot\xi(x)$ and $u\to u$,
 is not necessarily an automorphism of $R_{X,u}$, as it does not preserve the ideal $J$. However, using the $jet_0$ assumption, \S\ref{Sec.Notations.Rings.Germs},
   we can adjust it to an automorphism:  $\Phi_X: x\to x+t(f,u)^d\cdot c\cdot\xi(x)+t^2(f,u)^{d+1}$ and $u\to u$.
    Applying this $\Phi_X$ we get again equation \eqref{Eq.inside.proof.K.foliates.end.of.induction}.

 This completes the induction step in the general case.
\eei

\

Iterating this induction step one gets
\beq
\Phi_Y\Phi_X(f+tg)-f\in (x,u)^N\RpXu+Span_{R_{Y,u}}[(y,u)\cdot \{v_\bullet\}],u\big) \quad \text{ for } N\gg1.
\eeq
By the finite determinacy we conclude:
  $\Phi_Y\Phi_X(f+tg)-f\in  Span_{R_{Y,u}}[(y,u)\cdot \{v_\bullet\}],u\big)$.

 \item We have proved: $(f+tg,u)\in \cA\big(\{f\}+Span_{R_{Y,u}}[(y,u)\cdot \{v_\bullet\}],u\big)$ for $t$ an indeterminate.
   The key step was the solvability of the finite polynomial system \eqref{Eq.inside.proof.K.fibrations.system.of.IFT} in variables $\{c_{\um,\un}\}$.
    It defines a closed  algebraic subscheme in the (finite dimensional) affine space, $Z\sset Spec(\k[c_*])\times \k^1_t$.
     The projection of the algebraic germ  $(Z,o)\to (\k^1_t,o)$ is submersive.
      In particular, as $\k$ is infinite, the image of $Z$ is an infinite subset of $\k^1_t$.

Therefore $(f+tg,u)\in \cA\big(\{f\}+Span_{R_{Y,u}}[(y,u)\cdot \{v_\bullet\}],u\big)$ for an infinite set of values of $t\in \k^1_t$.

\item Inside the space $\MapX$ we study the intersection of the line $\{f+tg\}_{t\in \k}$ with the union of $\cA$-orbits,
 $\cA\big(\{f\}+Span_{R_{Y,u}}[(y,u)\cdot \{v_\bullet\}],u\big)$. By the
 finite determinacy we can pass to the finite jets, $jet_N(R_{X,u}):=\quots{R_{X,u}}{(u,x)^N}$. This is a finite-dimensional affine space.
  The transformations $\Phi_X,\Phi_Y$ of Step 2 were unipotent. Thus we have the algebraic, unipotent action of $\cA^{(1)}$ on the affine space
   $(\quots{R_{X,u}}{(u,x)^N})^{\oplus(p+r)}$.
   By Kostant-Rosenlicht theorem (see e.g. \cite[Theorem 2.11]{Ferrer Santos-Rittatore})
     its orbits are Zariski-closed. Moreover, $(f,u)+(Span_{R_Y}[(y)\cdot \{v_\bullet\}],0)$ is a linear subspace, and its $\cA^{(1)}$-orbit
    is Zariski closed as well.

By Step 3 the intersection of the line $\{f+tg\}_{t\in \k}$    with the orbit $\cA^{(1)}\big(\{f\}+Span_{R_{Y,u}}[(y,u)\cdot \{v_\bullet\}],u\big)$ is infinite.
 Therefore the whole line lies inside  $\cA^{(1)}\big(\{f\}+Span_{R_{Y,u}}[(y,u)\cdot \{v_\bullet\}],u\big)$. Namely,
  $(f+tg,u)\in \cA^{(1)}\big(\{f\}+Span_{R_{Y,u}}[(y,u)\cdot \{v_\bullet\}],u\big)$ for every $t\in \k^1$.
\epr
\eee

\section{Stable maps}\label{Sec.Stable.Maps}
Let $R_X$ be one of $\quots{\k[[x]]}{J}$, $\quots{\k\{x\}}{J}$, $\quots{\k\bl x\br}{J}$, see \S\ref{Sec.Notations.Rings.Germs}.

\subsection{Stability vs infinitesimal stability}\label{Sec.Stable.Maps.Stability.vs.Infinitesimal} Take a map  $f:X\stackrel{}{\to}(\k^p,o)$.
\bprop\label{Thm.Stable.vs.Infinitesimally.stable}
\bee
\item If $f$ is stable then it is infinitesimally  stable,  i.e. $T_\cA f=\RpX$,  i.e. $T^1_{\cA}f=0$.
\item
(If $char(\k)\!>\!0$\! then assume: $J=0$\! and $\k$\! is infinite.)
If $f$\! is infinitesimally stable then it is   stable.
\eee
\eprop
\bpr
\bee
\item
 Take a perturbation $v\in \RpX$ and the unfolding $f_t=f_o+t\cdot v$. As $f$ is stable, this unfolding is trivial. By theorem \ref{Thm.Unfolding.Triviality.Local}
  we get: $v\in T_{\cA_t}f_t$. The $t=0$ part of this condition is $v\in T_{\cA}f_o$. Therefore $T_{\cA}f_o=\RpX$.

\item  Take an unfolding $F(x,t)=(f_t(x),t)$ of $f= f_o$.
 As $T_\cA f_o=\RpX$ one gets $T_{\cA_t}f_o=\RpXt$. In particular the unfolding $F$ is   $\cA$-separable.
  Moreover, $f_o$ is $\cA_t$-finitely determined  as an element of $\RpXt$, see  \cite[\S7]{Kerner.Group.Orbits}.
 By lemma \ref{Thm.Unfoldings.Algebraic.Lemma} we get $T_{\cA_t}f_t=\RpXt$.

Therefore $\di_t f_t\in T_{\cA_t}f_t$. Now apply theorem \ref{Thm.Unfolding.Triviality.Local} to conclude: the unfolding $F$ is $\cA$-trivial.
\epr\eee

\bex\label{Ex.Stable.vs.Infinit.Stable}
For $R_X=\k\{x\}$ with $\k\in \R,\C$  and for $R_X=C^\infty(\R^n,o)$, this is the classical lemma \cite{Mather1968},
 see also Theorem 3.2 (pg. 62) of \cite{Mond-Nuno}.
\eex

\subsection{Stable maps are   unfoldings of their genotypes}\label{Sec.Stable.Maps.Unfoldings.of.Genotypes}
Take a map $F:\tX\to (\k^{p+r},o)$ of rank $r$.
  If $char(\k)>0$ then we assume: $\k$ is infinite and $J=0$ (i.e. the germ  $\tX$ is smooth).
Then (lemma \ref{Thm.prenormal.form.of.a.map}) $\tX\cong X\times(\k^r_t,o)$
 and $F(x,t)=(f(x)+c(x,t),t)$, for some $f\in (x)^2\cdot \RpX$ and $c(x,t)\in (x)\cdot (t)\cdot\RpXt$.
\bthe\label{Thm.Stable.Maps.unfoldings.of.their.genotypes}
 $F$ is stable iff $F$ is $\cA$-equivalent to the unfolding $(f+\sum^{r'}_{j=1} t_j v_j,t )$, where
  $t=(t_1,\dots,t_r)$, $r\ge r'$, $f\in (x)^2\cdot \RpX$
  is $\cK$-finite, and
   $\{v_1,\dots,v_{r'}\}\in (x)\cdot \RpX$ go to a set of generators of the $\k$-vector space $(x)\cdot T^1_\cK f$.
\ethe
In this case   the map $f$ is called ``the genotype" of $F$, see e.g. \cite[III.1.7]{AGLV}.

We emphasize: the generating tuple $\{v_\bullet\}$  can be chosen arbitrarily,  and is not necessarily minimal.
\\\bpr
 $\pmb\Lleftarrow$  By proposition \ref{Thm.Stable.vs.Infinitesimally.stable} it is enough to prove: $T_\cA F=R^{p+r}_{X,t}$.
  We have  $Span_{R_Y}\{v_\bullet\}+T_\cA f=\RpX$, cf. \eqref{Thm.TA.vs.TK}.

 Fix some generators $\{\xi_j\}$ of $Der_\k(R_X)$, then the generators of $Der_\k(R_{X,t})$ are  $\{\xi_j\}$ and $\{\di_{t_i}\}$.
The submodule $T_\cR F\sseteq R^{p+r}_{X,t}$ is generated by the matrix with columns $\{\xi_j F\},\{\di_{t_i}F\}$. We write these columns (as rows):
\beq
\{(\xi_j f+\sum t_i \xi_j  q_i,0,\dots,0)^t\}_{j=1,\dots,n},\quad \{q_{j-p}\cdot \he_1+\he_ j \}_{j=p+1,\dots,p+r'} ,\quad \he_{p+r'+1},\dots,\he_{p+r}.
\eeq
Here $\{\he_\bullet\}$ is the standard basis of $\k^{p+r}$.

By the direct check: $\RpX\oplus \zero\sseteq T_\cA F|_{t=0}$. And then $\zero\oplus R^{\oplus r}_X\sseteq T_\cA F|_{t=0}$. Altogether:
$T_\cA F|_{t=0}=R^{p+r}_X$.

By lemma \ref{Thm.Unfoldings.Algebraic.Lemma} we get $R^{p+r}_{X,t}=T_\cA F$.

\noindent$\pmb\Rrightarrow$
  By proposition \ref{Thm.Stable.vs.Infinitesimally.stable}, $F$ is infinitesimally-stable, i.e. $T_{\cA}F=R^{\oplus (p+r)}_{X,t}$.
Then $F$ is finitely $\cA$-determined, \cite[\S7]{Kerner.Group.Orbits}. Therefore  below we work over $\quots{\k[[x,t]]}{J}$,    with formal  $\cA$-transformations.

The map $F$ is an unfolding, thus we take $F$ in the pre-normal form $(f+\sum a_j(t,f)\cdot v_j,t)$, see Remark \ref{Rem.Pre-normal.form}.
 This unfolding is versal, therefore
   $Span_\k (\di_t \sum a_j(t,f)v_j|_{t=o})$ generates $(x)\cdot T^1_\cK f$. Thus
 $Span_\k (\di_t \sum a_j(t,o)v_j|_{t=o})$ generates $(x)\cdot T^1_\cK f$. And therefore the map $(\k^r,o)\ni t\to \{a_\bullet(t,o)\}\in (\k^{r'},o)$
  admits the (right) section. Thus, after an $\cR$-coordinate change on $t$-variables, we can assume $a_j(t,o)=t_j$.
  Combined with the $\cL$-transformation on the $(\k^r,o)$-part of the target, we get:
 $F\stackrel{\cA}{\sim}(f+\sum (t_j+(t)^d\cdot (f))v_j,t)$, for some $1\le d<\infty$.

The rest of the proof is induction.
 Apply  the $\cR$-transformation $ t_j+(t)^d\cdot (f)\rightsquigarrow t_j$ to get to $(f+\sum  t_j v_j,t+(t)^d\cdot (f))$.
  Now apply the $\cL$-transformation of type
  $(y,t)\rightsquigarrow (y,t-(t)^d\cdot(y) ))$ to get to  $(f+\sum  t_j v_j,t+(t)^{d+1}\cdot (x))$.
   Now apply the $\cR$-transformation $ t_j+(t)^{d+1}\cdot (x)\rightsquigarrow t_j$
   to get to $(f+\sum  t_j v_j+(t)^{d+1}\cdot(x),t)$. Iterate.
\epr

\bex
In the $\C$-analytic case,  $R_X=\C\{x\}$, this is Theorem 7.2 of \cite{Mond-Nuno}, for $C^\infty$-case see \cite[pg.200]{Martinet.1982}.
\eex
\bcor Any $\cK$-finite map $f\in \RpX$ admits a stable unfolding. Moreover, for a fixed number of parameters the stable
unfolding is unique up to $\cA$-equivalence.
\ecor
For $R_X=\C\{x\}$ this is Proposition 7.2 of \cite{Mond-Nuno}.

Given an ideal $I\sseteq (x)^2\sset R_X$ fix its generators, $I=R_X\{g_\bullet\}\sset R_X$. This defines the map $g:X\to (\k^p,o)$.
Define the Tjurina number of an ideal as $\tau(I):=codim_\cK g:=dim_k T^1_\cK g$. (This number does not depend on the choice of the generators.)
\bcor
 For every ideal $I\sseteq(x)^2\sset R_X$ with $\tau(I)<\infty$ there exists a stable germ $f\in \RpXt$ such that $(I,t)=(f)\sset R_{X,t}$.
\ecor
For $R_X=\C\{x\}$ this  Corollary 7.2 of \cite{Mond-Nuno}.

\subsection{Stable maps are determined by their local algebra}\label{Sec.Stable.Maps.Determined.by.local.Algebra}
\bprop (If $char(\k)>0$ then   assume: $J=0$  and $\k$ is an infinite field.)
 Two stable maps are $\cA$-equivalent iff they are $\cK$-equivalent.
\eprop
Recall, the stable maps are equivalent as maps, not as unfoldings.  See example ?? in \cite{Mond-Nuno}.
\bpr
(We prove the non-trivial direction, $\Lleftarrow$.) Stable maps are finitely determined, \cite[\S7]{Kerner.Group.Orbits},
 therefore we take $R_X=\quots{\k[[x]]}{J}$.

By theorem \ref{Thm.Stable.Maps.unfoldings.of.their.genotypes} we can take one map in the form  $F=(u,f +\sum u_j v_j )$.
  Here $f(x)\in (x)^2\cdot \RpX$ is $\cK$-finitely determined,    $u=(u_1,\dots,u_r)$, while $\{v_\bullet\}$ are
  sent to generators (not necessarily a basis) of the $\k$-vector space $(x)\cdot T^1_\cK f$.

By proposition \ref{Thm.Orbits.K.fibers.into.A} we can take the second map in the form
   $\tF=( u,f(x)+\sum  \tc_j(f,u) v_j)$, where $\tc_j\in (y,u)\cdot\k[[y,u]]$.
   As $\tF$ is infinitesimally stable, the map $u\to \tc(u,o)$ is invertible.
    Applying the $\cR_u$-transformation   $u_j+\tc_j(f,u)\rightsquigarrow u$ we get: $\tF\stackrel{\cR}{\sim} ( u+c(f,u) ,f(x)+\sum   u_j  v_j(x))$.
     Then by $\cL$-transformation we get to
  $\tF\sim ( u+q(u,x) ,f(x)+\sum   u_j  v_j(x))$, where $q(u,x)\in (u)\cdot(x)\cdot R^{\oplus r}_{x,u}$. Now again an $\cR$
   transformation on $u$ gives  $\tF\sim ( u ,f(x)+\sum   u_j  v_j(x)+(u)\cdot (x))$.
 Therefore $\tF$ is an unfolding of $f$. Finally, bring $\tF$ to the pre-normal form, theorem \ref{Thm.Stable.Maps.unfoldings.of.their.genotypes}.
\epr

\section{Results of Mather-Yau/Scherk/Gaffney-Hauser type}\label{Sec.M.Y.G.H.}
 It is well known that the mapping $f:(\k^n,o)\to (\k^p,o)$ is determined (up to $\cR,\cK,\cA$-equivalence)
 by its ``behaviour" at the critical/singular/instability locus. One  way to obtain a precise statement was given in \cite{Mather-Yau} for $\cK$-equivalence
  of functions ($p=1$),   in \cite{Scherk} for $\cR$-equivalence of functions,
  and extended in \cite{Gaffney-Hauser} to $\cK$ and $\cA$-equivalences of maps ($p\ge1$).
   The initial proofs were $\C$-analytic. The first extension  to zero/positive characteristic was
   done for $p=1$, $\cR,\cK$ in \cite{Greuel-Pham.2019} (in the case of isolate critical point).
    We establish the general statements/strengthen the known results.

\subsection{The $\cK$-version in zero characteristic}\label{Sec.M.Y.G.H.char.0}
Consider modules over rings, $M_j\in mod\text{-}R_j$, with the following notion of isomorphism:
 $M_1\sim M_2$ if $M_1\cong \phi^* M_2$ for an isomorphism of rings $\phi: R_1\isom{} R_2$.

Below we consider the module $T^1_\cK f$ up such isomorphisms.
  For $p=1$ the data ``$  T^1_\cK f $ up to isomorphism" is equivalent to the data of $\k$-algebra $T^1_\cK f=\quots{R_X}{(f)+Jac( f)}$.

Let $R_X$ be one of $\quots{\k[[x]]}{J}$, $\quots{\k\{x\}}{J}$, $\quots{\k\bl x\br}{J}$, with $char(\k)=0$, $\k=\bar\k$.
 Take a map $f:X\to (\k^p,o)$.
\bthe
\bee
\item  The $\cK$-type of $f$ is determined by the isomorphism type of  $T^1_{\cK^{(0)}} f$.
\item If  $f$ is $\cK$-finite then the $\cK$-type of $f$ is determined by the isomorphism type of  $T^1_{\cK} f$.
\eee
\ethe
\noindent More precisely,  if
$T^1_{\cK^{(0)}} f_0\sim T^1_{\cK^{(0)}} f_1$ (in the sense as above), then $f_0\stackrel{\cK^{(0)}}{\sim}f_1$.
 (And similarly for part 2.)
\\\bpr The proofs of parts 1,2 are the same except for Step 3.i. To simplify notations we  work mostly with $T_\cK$ and $T^1_{\cK} f$.

\bee[{\bf Step 1.}]
\item
\bee[\bf i.]
\item It is enough to establish the formal case, $R_X=\quots{\k[[x]]}{J}$.
 For the statements over $\quots{\k\{x\}}{J}$, $\quots{\k\bl x\br}{J}$ one applies Artin approximation.

\item Given an isomorphism $T^1_{\cK } f_0\cong \phi^* T^1_{\cK } f_1$ coming from an automorphism $\phi\circlearrowright R_X$,
 replace $f_1$ by $f_1\circ \phi^{-1}$ to
 get an isomorphism of $R_X$-modules $T^1_{\cK } f_0\cong T^1_{\cK } f_1$. After an $R_X$-linear automorphism of the module $\RpX$
  we can assume: $T_\cK f_0= T_\cK f_1\sseteq \RpX$.

Similarly, in the first part we can assume:  $T_{\cK^{(0)}} f_0= T_{\cK^{(0)}} f_1\sset  \RpX$.

\eee
 \item Extend the ring, $R_{X,t}:=R_X[t]$, note that  $R_{X,t}$ is non-local.  Take the unfolding $f_t:=f_0+t(f_1-f_0)$.
  We have the tangent spaces $T_{\cK_t}$, $T_{\cK_t}f_t$ and $T_{\cK^{(0)}_t}f_t$,   see \S\ref{Sec.Triviality.of.Unfolding.Nonlocal.Base}.

 We claim: $\di_t f_t|_{t_o}\in T_{\cK _t} f_t|_{t_o}$ for $t_o\in \k^1\smin \{finite\ set\}$.
  (For the first part: $\di_t f_t|_{t_o}\in T_{\cK^{(0)}_t} f_t|_{t_o}$ for $t_o\in \k^1\smin \{finite\ set\}$.)
  More precisely, we claim: $\di_t f_t \in T_{\cK _t} f_t $ over the factor ring $\k[t][g^{-1}]$, for some polynomial $0\neq g\in \k[t]$.

 Obviously $\di_t f_t=f_1-f_0\in  T_{\cK } f_1+ T_{\cK _t} f_0= T_{\cK } f_0$. Thus it is enough to prove:
 $T_{\cK _t} f_0= T_{\cK _t}f_{t_o}$ for all $t_o\in \k^1$ except for a finite set.

\bee[\bf i.]
\item   We have: $T_{\cK _t} f_t\sseteq T_{\cK _t}f_0+t(T_{\cK _t} f_0+T_{\cK _t} f_1)=T_{\cK _t} f_0$.
\item Similarly: $T_{\cK _t} f_0\sseteq T_{\cK _t}f_t+t\cdot T_{\cK _t}(f_1-f_0) \sseteq   T_{\cK _t}f_t+t\cdot T_{\cK _t} f_0$.
 These $R_{X,t}$-modules
 are finitely generated. Localize at the ideal $(x,t)$ and apply Nakayama over the local ring  $(R_{X,t})_{(x,t)}$.
 We get: $(T_{\cK _t} f_0)_{(x,t)}\sseteq  (T_{\cK _t} f_t)_{(x,t)}$.
\item Similarly one has: $T_{\cK _t} f_1\sseteq T_{\cK _t}f_t+(1-t)\cdot  T_{\cK _t}(f_1-f_0)
   \sseteq T_{\cK _t}f_t+(1-t)\cdot T_{\cK _t} f_1$.
 Therefore $(T_{\cK _t} f_0)_{(x,1-t)}=(T_{\cK _t} f_1)_{(x,1-t)}\sseteq  (T_{\cK _t} f_t)_{(x,1-t)}$.
\item Take the quotient module
    $M:=\quots{T_{\cK _t}f_0}{T_{\cK _t}f_t}$. Thus $M$ is finitely generated over $R_X[t]$.
     Its localizations at two points vanish,  $M_{(t,x)}=0$ and $M_{(1-t,x)}=0$.
    Therefore its support,  $Supp(M)\sseteq Spec(R_X[t])$, does not contain the points $V(t,x),V(1-t,x)\in Spec(R_X[t])$. But $Supp(M)$ is
        a Zariski-closed subset. Therefore $Supp(M)\cap V(x)\sset \k^1_t$ is a finite set of points.
\eee
Altogether, we have proved: $\di_t f_t|_{t_o}\!\in\! T_{\cK _t} f_t|_{t_o}$ (resp. $\di_t f_t|_{t_o}\!\in\! T_{\cK^{(0)}_t} f_t|_{t_o}$)
  for $t_o\!\in\! \k^1\!\smin\! \{finite\ set\}$.

\item
\bee[\bf i.]
\item ({\bf The case $\k=\C$.}) We have the unfolding $f_t$ over $\C$, and
  $\di_t f_t|_{t_o}\in T_{\cK _t} f_t|_{t_o}$ for $t_o\in \k^1\smin \{finite\ set\}$.
Take a path $\ga:0\stackrel{}{\rightsquigarrow}1$ in $\C$ that avoids this finite set.
 By theorem \ref{Thm.Unfolding.Triviality.Local} the unfolding $f_t$ is locally $\cK$ (resp. $\cK^{(0)}$)-trivial at each point of the path.

The $\cK^{(0)}$-trivialization (in the first case) preserves the origin of $X$ by definition.
The $\cK $-trivialization (in the second case) preserves the origin of $X$ as $V(f_0),V(f_1)$ have isolated singularities.

As the path is compact we take a finite cover $\ga=\cup_i U_i$, such that $f_t$ is trivial on each $U_i$.
And  by the connectedness of the
 path we get: $f_0\stackrel{\cK }{\sim}  f_1 $, resp. $f_0\stackrel{\cK^{(0)}}{\sim}  f_1 $.

\item ({\bf The general case, $\k=\bar\k$, $char(\k)=0$.}) Pass to the finite jets, $jet_d R_X:=\quots{R_X}{(x)^{d+1}}$.
 Then for  $t\in \k^1\smin \{finite\ set\}$ we have:
   $\di_t jet_d f_t=jet_d\di_t  f_t\in jet_d T_{\cK _t} f_t\sseteq  T_{\cK_t(jet_d R_X)} jet_d f_t$. Explicitly, we have
    $\di_t jet_d f_t=jet_d \xi_X(jet_d f_t)+jet_d(U) \cdot jet_d f_t$.
 This holds for almost all $t\in \k^1$, i.e. over $\k[t][g^{-1}]$.

Taylor-expand the elements $jet_d \xi_X, jet_d U , jet_d(J)$ up to order $d$. Let $\{C_\bullet\}$ be the (finite) set of the coefficients.
  Then $\Q\{C_\bullet\}$ is a finite field extension of $\Q$. Therefore this extension can be (re-)embedded,
   $\ep:\Q\{C_\bullet\}\stackrel{ }{\hookrightarrow} \C$.

Over $\C$ we still have:  $\di_t jet_d \ep(f_t)\in  T_{\cK _t(\ep(jet_d R_X))} jet_d \ep(f_t)$  for almost all $t\in \k^1$.
 By part i. we get equivalence over $\C$: $jet_d \ep(f_0)\stackrel{\cK_t(\ep(jet_d R_X))}{\sim}jet_d \ep(f_1)$.
  (resp. $jet_d \ep(f_0)\stackrel{\cK^{(0)}_t(\ep(jet_d R_X))}{\sim}jet_d \ep(f_1)$).
 Then   Lemma \ref{Thm.Equivalence.Change.base.field} gives equivalence over $\k$:
     $jet_d  f_0 \stackrel{\cK_t( jet_d R_X )}{\sim}jet_d  f_1 $.

This holds for each $d\gg1$. Therefore the condition $f_1\in \cK _t(f_0)$, which is an implicit function equation, has an order-by-order solution.
 By Pfister-Popescu theorem, \cite{Pfister-Popescu}, we get the   equivalence  $  f_0 \stackrel{ }{\sim}  f_1 $ over
  $\quots{\k[[x]]}{J}$.\epr
\eee
\eee
\beR
\bee[\bf i.]
\item The statement does not hold if $\k$ is not algebraically closed. E.g. suppose $\sqrt[d]{a}\not\in \k$, for some $d\ge4$,
 and compare $x^d_1+x^d_2$ to $x^d_1+a\cdot x^d_2$.
\item If the singularity is not ``of isolated type" then the $\cK$-type of $f$ is not determined by the module $T^1_\cK f$. See \cite[\S4]{Gaffney-Hauser}.
 Recall the standard example. Let $f_t(x,y)=xy(x-y)(x-(z-t)y)$. Here $T_\cK f_t=(x^3,x^2y,xy^2)$, independent of $t$.
  This family is trivialized by  the coordinate change  $\phi:(x,y,z)\to(x,y,z+t)$, which does not preserve the origin.
   Thus $\phi$ is not an automorphism of the local ring $\k[[x,y,z]]$. And no automorphism of $\k[[x,y,z]]$ can trivialize this family,
    as any automorphism will preserve the cross-ratio of the four planes of $V(f_t)\sset (\k^3,o)$.
\item Instead of the module $T^1_\cK f$ we could take the $R_X$-module $T^1_\cR f$, getting the similar statement.
\eee
\eeR

\subsection{The $\cA$-version in zero characteristic} \
\\\parbox{14cm}{Below we consider $T^1_\cA f$ as a mixed $(R_Y,R_X)$-module. We write $T^1_\cA f\sim T^1_\cA \tf$ if there exists an isomorphism
of algebras $(\Phi_X,\Phi_Y)$ (see the diagram)
 that induces the isomorphism of $R_Y$-modules
 $\Phi_Y': T^1_\cA f\isom{} T^1_\cA \tf$.
}  \quad\quad\quad$ \bM R_X  \isom{\Phi_X} R_X\\\uparrow f^\#\ \tf^\#\uparrow\\R_Y \underset{}{\isom{\Phi_Y}} R_Y\eM $

\vspace{-0.5cm}  \bex
\bee[\bf i.]
\item
 If $f\stackrel{\cA}{\sim}\tf$, i.e. $\tf=\Phi_Y\circ f\circ\Phi_X$, then $T^1_\cA f\sim T^1_\cA \tf$.
  Indeed, $T^1_\cA(\tf\circ\Phi_X^{-1})=\quots{\RpX}{T_\cR(\Phi_Y f)+T_\cL(\Phi_Y f)}=\quots{\RpX}{\Phi'_Y T_\cR f+T_\cL  f}\isom{\Phi_Y'}T^1_\cA f.$
  \item Suppose $T^1_\cA f\sim T^1_\cA \tf$. Using the corresponding morphism of algebras we get $T^1_\cA(\Phi_Y\circ f\circ\Phi_X)= T^1_\cA \tf$, i.e.
 $T_\cA(\Phi_Y\circ f\circ\Phi_X)= T_\cA \tf\sset \RpX$.
\eee \eex

\bthe Let $R_X=\quots{\k[[x]]}{J}$, with $\k=\bar\k$ and $char(\k)=0$. Suppose $\sqrt{I}=\cm$.
\bee
\item The $\cA$-type of $f$ is determined by the mixed module type of  $T^1_{\cA^{(0)}} f$.
\item If $f$ is $\cA$-finite\! and\! not\! stable\! then\! the $\cA$-type of $f$ is determined by the mixed module type of  $T^1_{\cA } f$.
\eee
\ethe
More precisely, if $T^1_{\cA^{(0)}} f_0\sim T^1_{\cA^{(0)}} f_1$, then $f_0\stackrel{\cA^{(0)}}{\sim}f_1$. (And similarly in case 2.)
\\\bpr By the example above we can assume $T_{\cA^{(0)}} f_0=T_{\cA^{(0)}} f_1\sset\RpX$ (resp. $T_{\cA} f_0=T_{\cA} f_1\sset\RpX$).
 The proof is similar to the $\cK$-case and is inductive.
 We apply a sequence of transformations $g_d\in \cA^{(0)}(jet_d(R_X))$ and verify
 $jet_d f\stackrel{g_d\cdots g_1}{\sim}jet_d \tf$. Moreover, the transformation can be chosen to satisfy $jet_{d-1}(g_d)=Id$.

 Fix any $d\ge1$ and replace $R_X$ by the Artinian ring $\quots{R_X}{(x)^{d+1}}$. Assume $jet_{d-1}(f_0)=jet_{d-1}(f_1)$.

\bee[\bf Step 1.]
\item As in the $\cK$-case we define  $f_t:=f_0+t(f_1-f_0)$ and $R_{X,t}:=R_X\otimes\k[t]$, $T_{\cA_t}:=T_\cA\otimes\k[t]$.
 We claim: $\di_t f_t\in T_{\cA^{(j_d)}_t}f_t$ for $t\in \k^1\smin $(finite set), i.e. this holds over the fraction ring $\k[t][S^{-1}]$.
 Here $j_d\ge-1$ and $\lim_{d\to \infty}j_d=\infty$.

We have: $\di_t f_t=f_1-f_0\in (T_{\cA^{(0)}} f_0+T_{\cA^{(0)}} f_1)\cap \cm^d\cdot \RpX=T_{\cA^{(0)}} f_0\cap \cm^d\cdot \RpX$.
 (For part two one gets: $\di_t f_t\in T_{\cA} f_0\cap \cm^d\cdot \RpX$.)

By the ``mixed Artin-Rees" lemma, \cite[\S7]{Kerner.Group.Orbits}, one has  $T_{\cA} f_0\cap \cm^d\cdot \RpX\sset T_{\cA^{(j_d)}} f_0$, where $\lim_d j_d=\infty$.

Therefore, as in the $\cK$-case, it is enough to verify $T_{\cA^{(j_d)}_t} f_t= T_{\cA^{(j_d)}_t}f_0$ for $t\in \k^1\smin $(finite set).
 \bei
 \item
We have $T_{\cR^{(j_d)}_t}f_t\sseteq (1-t)\cdot T_{\cR^{(j_d)}_t}f_0+t\cdot  T_{\cR^{(j_d)}_t}f_1\sseteq  T_{\cA^{(j_d)}_t}f_0$. In addition
  \beq
  T_{\cL_t}f_t\sseteq \{\sum_{i,j\ge0} q_i(f_0)\cdot v_j(f_1)\}\sseteq  {  T_{\cL_t}f_0+  T_{\cL_t}f_1}\sseteq
      {  T_{\cA_t}f_0}.
  \eeq
(Note that $R_X$ is Artinian now.) And similarly $T_{\cL^{(j_d)}_t}f_t\sseteq T_{\cA^{(j_d)}_t}f_0$.

\item As in the $\cK$-case we have: $T_{\cA^{(j_d)}_t} f_0\sseteq T_{\cA^{(j_d)}_t} f_t+t\cdot T_{\cA^{(j_d)}_t} f_0$.
 These modules are finitely generated over $R_{Y,t}$, because $R_X$is Artinian.
  Localize at $(y,t) $ to get:
$(T_{\cA^{(j_d)}_t} f_0)_{(t,y)}\sseteq (T_{\cA^{(j_d)}_t} f_t)_{(t,y)}$.

Similarly one gets  $(T_{\cA^{(j_d)}_t} f_1)_{(t,y)}\sseteq (T_{\cA^{(j_d)}_t} f_t)_{(t,y)}$.

\item
  As in the $\cK$-case take the quotient $M:=\quots{T_{\cA^{(j_d)}_t} f_0}{T_{\cA^{(j_d)}_t} f_t}$. This is an $R_{Y,t}$-module.
   Thus: $ M_{(t,y )}=0, $ and $M_{(1-t,y )}=0$. Hence $Supp(M)\sset V(y)\sset Spec(R_{Y,t})$ is a finite subset.

\eei
\item
\bee[\bf i.]
\item  ({\bf The case $\pmb{ \k=\C}$.})  As in the $\cK$-case take  a path inside $\C_t\cong V(y)\sset (\C^p_y,o)\times(\C^1_t,o)$ avoiding the finite subset $Supp(M)$. Along this
 path one has $\di_t f_t\in T_{\cA^{(j_d)}_t} f_t$. Now use part $\cA$ of theorem \ref{Thm.Unfolding.Triviality.Local}
   to get $f_0\stackrel{\cA^{(j_d)}}{\sim}f_1$.
  Note that the ring $R_X$ is Artinian, therefore the formal and ordinary equivalence coincide.

For part 2 we remark that the $\cA$-trivialization preserves the origins of $(\k^n,o)$, $(\k^p,o)$, as $f_0,f_1$ are $\cA$-finite, i.e. their
 instability locus is just one point.

\item ({\bf The general case, $\pmb{ \k=\bar\k, char(\k)=0}$.}) By the same Lefschetz-type arguments (and lemma \ref{Thm.Equivalence.Change.base.field})
  we get $f_0\stackrel{\cA^{(j_d)}}{\sim}f_1$.
 \eee
\eee
Now combine these constructed transformations to get the (formal) group element $g:=\lim_d(g_d\cdots g_1)\in \cA^{(0)} $
 (resp. $g:=\lim_d(g_d\cdots g_1)\in \cA $), satisfying $g \tf=f$.
\epr

\subsection{The case of arbitrary characteristic}\label{Sec.M.Y.G.H.any.char}
 Let $R_X$ be one of $\quots{\k[[x]]}{J}$, $\quots{\k\{x\}}{J}$, $\quots{\k\bl x\br}{J}$, here $\k$ is any  field.
  If $char(\k)>0$ then we assume the $jet_0$ condition of \S\ref{Sec.Notations.Rings.Germs}.
 Take a map $f:(\k^n,o)\to (\k^p,o)$, thus $0\neq f\in \cm\cdot \RpX$.
 The order of $f$
  is the largest $ord(f)\in \N$ satisfying $(f)\sseteq \cm^{ord(f)}$. If $ord(f)\le 2$ then we take $\cm^{ord(f)-2}=R_X$.
\bthe
\bee[\bf 1.]
 \item ($\cK$-case) Suppose an ideal $\ca\sseteq \cm^2$ satisfies:  $\ca^2\cdot\cm^{ord(f)-2}\cdot \RpX \sseteq \cm\cdot \ca\cdot T_\cR f+\cm\cdot (f)\cdot\RpX$.
 Then the $\cK$-type of $f$ is determined by the   $\k$-algebra   $\quots{R_X}{(f)+\ca\cdot \ca_\cR}$.
 \item ($\cA$-case, $\k$-infinite) Suppose an ideal $\ca\sseteq \cm^2$ satisfies:
  $\ca^2\cdot\cm^{ord(f)-2}\cdot\RpX \sseteq   \ca\cdot T_\cR f+f^{\#}(y)^2\cdot T_\cL f$.
 Then the $\cA$-type of $f$ is determined by the $\k$-algebra   $\quots{R_X}{f^{\#}(y)+\ca^2}$.

\item   ($\cR$-case) Suppose an ideal $\ca\sseteq \cm$ satisfies:  $\ca^2\cdot \cm^{ord(f)-2}\cdot\RpX\sseteq \cm\cdot \ca\cdot T_\cR f$.
 Then the $\cR$-type of $f$ is determined by the $\k[f]$-algebra   $\quots{R_X}{\ca\cdot \ca_\cR}$.
 \eee
\ethe
Here is the explicit form of the statement.  Fix some maps $f,\tf\in \cm\cdot \RpX$ and ideals $\ca^f,\ca^\tf\sseteq\cm^2$.
\bei
\item ($\pmb\cK$) Suppose
$(\ca^f)^2 \cdot \cm^{ord(f)-2}\cdot\RpX \sseteq \cm\cdot \ca^f\cdot T_\cR f+\cm\cdot (f)\cdot\RpX$ and
$(\ca^\tf)^2 \cdot\cm^{ord(\tf)-2}\cdot\RpX \sseteq \cm\cdot \ca^\tf\cdot T_\cR \tf+\cm\cdot (\tf)\cdot\RpX$.
If the $\k$-algebras are isomorphic,   $\quots{R_X}{(f)+\ca^f\cdot \ca_\cR^f}\cong    \quots{R_X}{(\tf)+\ca^\tf\cdot \ca^\tf_\cR}$,
 then $f\stackrel{\cK}{\sim}\tf$.

\item ($\pmb\cA$) Suppose
$(\ca^f)^2 \cdot \cm^{ord(f)-2}\cdot\RpX \sseteq \cm\cdot \ca^f\cdot T_\cR f+f^{\#}(y)^2\cdot T_\cL f$ and
$(\ca^\tf)^2 \cdot \cm^{ord(\tf)-2}\cdot\RpX \sseteq \cm\cdot \ca^\tf\cdot T_\cR \tf+\tf^{\#}(y)^2\cdot T_\cL \tf$.
If the $\k$-algebras are isomorphic,   $\quots{R_X}{f^{\#}(y)+(\ca^f)^2}\cong    \quots{R_X}{\tf^{\#}(y)+(\ca^\tf)^2}$, then $f\stackrel{\cA}{\sim}\tf$.

  \item ($\pmb\cR$)
Suppose $(\ca^f)^2 \cdot \cm^{ord(f)-2}\cdot\RpX\sseteq \cm\cdot \ca^f\cdot T_\cR f$ and
$(\ca^\tf)^2\cdot \cm^{ord(\tf)-2}\cdot\RpX\sseteq \cm\cdot \ca^\tf\cdot T_\cR \tf$. Suppose the $\k[f]$-algebra $\quots{R_X}{ \ca^f\cdot \ca_\cR^f}$
 is isomorphic to the $\k[\tf]$-algebra $\quots{R_X}{ \ca^\tf\cdot \ca^\tf_\cR}$. Then $f\stackrel{\cR}{\sim}\tf$.
\eei
We do not assume isolated singularities or that the germ $V(f)\sset X$ is a complete intersection.
\\\bpr
\bee[\bf 1.]
\item Given an isomorphism of $\k$-algebras $\quots{R_X}{(f)+\ca^f\cdot \ca_\cR^f}\isom{} \quots{R_X}{(\tf)+\ca^\tf\cdot \ca^\tf_\cR }$
 take its representative $\phi:R_X\stackrel{}{\to}R_X$.
 Thus $\phi$ is invertible (hence an isomorphism of  $\k$-algebras) and sends $(f)+\ca^f\cdot \ca_\cR^f$ to $(\tf)+\ca^\tf\cdot \ca_\cR(\tf)$.
  Therefore after a coordinate change ($\cR$-equivalence) we can assume:
  \beq\label{Eq.inside.proof.K.case}
  (f)+\ca^f\cdot \ca_\cR^f=(\tf)+\ca^\tf\cdot \ca_\cR(\tf)\sset R_X.
\eeq
Note that the assumption $\ca^2(f)\cdot\cm^{ord(f)-2}\cdot\RpX \sseteq \cm\cdot \ca^f\cdot T_\cR f+\cm(f)\cdot \RpX $ is preserved under $\cK$-equivalence,
 see \cite[\S5]{Kerner.Group.Orbits}.

By the initial assumption $\ca^f,\ca^\tf\sseteq\cm^2$.
  Then equation \eqref{Eq.inside.proof.K.case} gives: $rank(f)=rank(\tf)$. By $\cK$-equivalence we can assume $f=(x_1,\dots,x_r,\cf)$,
  with $(\cf)\sseteq (x_{r+1},\dots,x_n)$. (And similarly for $\tf$.) Therefore the whole question is reduced to the case $(f),(\tf)\sseteq\cm^2$.

Now apply $GL(p,R_X)$ transformations to $f$ and $\tf$ to get:
 $f- \tf\in  (\cm\cdot (f)+\ca^f\cdot \ca^f_\cR )\cdot \RpX$. Finally, by  \cite[\S5]{Kerner.Group.Orbits},
   we have $\cK f\supseteq \{f\}+(\cm\cdot (f)+\ca^f\cdot \ca^f_\cR )\cdot \RpX$.      In particular, $\tf\in \cK f$.

\item
As in the $\cK$-case we lift the isomorphism $\quots{R_X}{f^{\#}(y)+\ca^2(f)}\isom{}    \quots{R_X}{\tf^{\#}(y)+\ca^2(\tf)}$ to an isomorphism
 of $\k$-algebras
$\phi:R_X\stackrel{}{\isom{}}R_X$ sending $f^{\#}(y)+ (\ca^f)^2$ to $\tf^{\#}(y)+(\ca^\tf)^2 $. Then by a coordinate change ($\cR$) we can assume
 $f^{\#}(y)+(\ca^f)^2=\tf^{\#}(y)+(\ca^\tf)^2\sset R_X$. Therefore $\tf^{\#}(y)\in f^{\#}(y)+(\ca^f)^2$.
  Applying $GL(p,\k)$ transformations to $\tf$   we can assume:
 $f-\tf\in f^{\#}(y)^2T_\cL f+(\ca^f)^2\cdot \RpX$. Finally, by \cite[\S7]{Kerner.Group.Orbits},
   we have $\cA \tf\supseteq \{ f\}+  f^{\#}(y)^2\cdot T_\cL f+(\ca^f)^2\cdot \RpX$.      In particular, $f\in \cA g$.

\item
We have the isomorphism  $\quots{R_X}{ \ca^f\cdot \ca_\cR^f}\cong \quots{R_X}{ \ca^\tf\cdot \ca^\tf_\cR}$, compatible with the isomorphism
  $\k[f]\cong\k[\tf]$. The later isomorphism (after an automorphism of $\k[\tf]$) can be taken as $f_i\to \tf_i$.
   Then  $\quots{R_X}{ \ca^f\cdot \ca_\cR^f}\cong \quots{R_X}{ \ca^\tf\cdot \ca_\cR(\tf)}$ becomes an isomorphism of $\k[f]$-algebras.

Lift it to an isomorphism $R_X\stackrel{\phi}{\to}R_X$.
  Therefore by a coordinate change ($\cR$-equivalence) we can assume:  $\ca^f\cdot \ca_\cR^f=\ca^\tf\cdot \ca_\cR(\tf)\sset R_X$. Moreover, we have:
   $\tf_i=\tf_i\cdot 1\in f_i+\ca^f\cdot \ca^f_\cR\cdot \RpX $.
 By \cite[\S4]{Kerner.Group.Orbits} we get $\cR f\supseteq \{\tf\}+\ca^\tf\cdot \ca^\tf_\cR\cdot \RpX$.
     In particular, $\tf\in \cR f$.
\epr
\eee

\bex  Take $p=1$ and assume $J=0$, i.e. $X\cong(\k^n,o)$.
\bee[\bf i.]
\item (The $\cK$-case.)
 Suppose $Jac(f)\cdot \cm\cdot \ca+\cm\cdot (f)\supseteq \ca^2\cdot \cm^{ ord(f)-2}$.
  Then the $\cK$-type of $f$ is determined by the $\k$-algebra $\quots{R_X}{(f)+\ca\cdot Jac(f)}$. The singularity can be non-isolated,
     we do not assume $\sqrt{\ca}=\cm$.

As a particular case take $\ca=\cm^d$ with $d\ge2$.
Suppose $Jac(f)\cdot \cm^{d+1}+\cm\cdot (f)\supseteq \cm^{2d+ord(f)-2}$. Then the $\cK$-type of $f$ is determined by  $\quots{R_X}{(f)+\cm^d\cdot Jac(f)}$.
 Compare this to \cite[Theorem 2.2]{Greuel-Pham.mather-yau.char.positive} (for $R_X=\k[[x]]$):

 {\em If $Jac(f)\cdot \cm^2+\cm\cdot (f)\supseteq \cm^{\frac{d}{2}+ord(f)}$ then $\cK f$ is determined by the $\k$-algebra $\quots{R_X}{(f)+\cm^d\cdot Jac(f)}$.}
\\
Their assumption implies (is stronger than) the condition $Jac(f)\cdot \cm^\frac{3d}{2}+\cm^{\frac{3d}{2}-1}\cdot (f)\supseteq \cm^{2d+ord(f)}$,
 which is much stronger than ours.
\item (The $\cR$-case.)
Suppose $Jac(f)\cdot \cm\cdot \ca \supseteq  \ca^2\cdot \cm^{ ord(f)-2}$. Then the $\cR$-type of $f$ is determined by the $\k[f]$-algebra
  $\quots{R_X}{ \ca\cdot Jac(f)}$. Again,      we do not assume $\sqrt{\ca}=\cm$.

As a particular case take take $\ca=\cm^d$  with $d\ge2$.
Suppose $Jac(f)\cdot \cm^{d+1} \supseteq \cm^{2d+ord(f)-2}$. Then the $\cR$-type of $f$ is determined by  $\quots{R_X}{ \cm^d\cdot Jac(f)}$.
 This strengthens (and extends)   \cite[Theorem 2.4]{Greuel-Pham.mather-yau.char.positive} (for $R_X=\k[[x]]$).
 \eee
\eex
The $\cR,\cK$-cases for $p>1$ and the $\cA$-case are new.
\beR
\bee[\bf i.]
\item
One would like a stronger statement, e.g. (for $p=1$) of the form ``the $\cK$-type is determined by the $\k$-algebra $\quots{R}{(f)+\cm^d\cdot Jac(f)}$,
 where $d$ depends on $\k$ and $dim(R_X)$, but not on $f$". This is impossible due to the following example:
 \beq
 char(\k)=p,  R_X=\k[[x,y]]     \text{ and }   f(x,y)=x^{p+1}+y^{pN+1},  \tf(x,y)=f(x,y)+x^py^{pd},\quad \text{ for } \quad N>pd.
 \eeq
Here $Jac(f)=Jac(\tf)$ and $\cm^j\cdot Jac(\tf)\supset \cm^d\cdot x^p\ni x^py^{pd}$ for any $j\le d$. Therefore $(f)+\cm^j\cdot Jac(\tf)=(\tf)+\cm^j\cdot Jac(g)$ for $j\le d$.
 But $f\not\stackrel{\cK}{\sim}g$, e.g. because the monomial $x^py^{pd}$ lies under the Newton diagram of $f$.

\item
(For $p=1$) Recall that for $R=\C\{x\}$ the  $\cR$-type of $f$ is not determined just by the $\k$-algebra structure of
 $\quots{R}{Jac(f)}$. Moreover, the $\cR$-type is not determined by the $\k$-algebra  $\quots{R}{Jac(f)^d}$, for any $d$.
 Indeed, suppose $f\stackrel{\cR}{\not\sim}\tf$, but  $f\stackrel{\cK}{\sim}\tf$ via $\tf(x)=c\cdot f(\phi(x))$, with $c\in \k$.
  (An explicit example is $f_t(x_1,x_2)=x^4_2+x^5_1+t\cdot x^2_2x^3_1$,   see   \cite[pg.133]{Gr.Lo.Sh}.) Then $Jac(f_t)^d=Jac(f_0)^d$ for all $d\ge1$.
\eee
\eeR

\appendix
\section{Separability of unfoldings}\label{Sec.Separable.Unfoldings}
 Below we assume $\k=\bar\k$.
To an unfolding $F=(f_t(x,u),t)$ we associate the element $[f_t]\in T^1_{\cG_t}f_o$. Define the group action
\beq
G:=Aut_{(\k^r,o)}:=Aut_\k(\K)\circlearrowright\RpXt  \quad \text{ by } \quad  q(x,t)\to q(x,\phi(t)).
\eeq
  This action preserves the embedding $T_{\cG_t}f_o\sseteq \RpXt$ and hence descends to the action $G\circlearrowright T^1_{\cG_t}f_o$.
 We get the group orbit $[Gf_t]\sseteq T^1_{\cG_t}f_o$ and the orbit map $G\stackrel{f_t}{\to} T^1_{\cG_t}f_o$ by $g\to [g\cdot f_t]$.

 Take an ideal $\cb\sset R_{X,t}$ satisfying $\sqrt{\cb}=(x,t)$ and pass to the finite jets, $jet_\cb(R_{X,t}):=\quots{R_{X,t}}{\cb}$.
  This is a finite-dimensional $\k$-vector space. Similarly one has $jet_\cb(\RpXt)$, $jet_\cb(T^1_{\cG_t}f_o):=\quots{\RpXt}{T_\cG f_o+\cb\cdot \RpXt}$
   and $jet_\cb(G):=Aut_\k(\quots{\k[[t]]}{\k[[t]]\cap \cb})$. We get the (regular) action of the affine algebraic group on the affine space,
    $jet_\cb(G)\circlearrowright jet_\cb(T^1_\cG f_o)$.

\bed (Suppose $\k=\bar\k$.)
The orbit map $G\stackrel{f_t}{\to}T^1_{\cG_t}f_o$ is called separable if all its finite jets,
 $\{jet_\cb(G)\stackrel{jet_\cb(f_t)}{\to}jet_\cb T^1_\cG f_o\}_\cb$, are separable as morphisms onto their images.
\eed
(The morphism $jet_\cb(f_t)$ is called separable if the corresponding field extensions are separable, see e.g. \cite{Ferrer Santos-Rittatore}.)

Take the tangent space, $T_G:=(t)\cdot Der_{(\k^r,o)}$, and its finite jets, $\{jet_\cb T_G\}_\cb$. Take the map from the jet(image tangent space) to
 the tangent space of the jet-orbit, $jet_\cb(T_G f_t)\to T jet_\cb(G f_t)$. We write $T_G f_t\twoheadrightarrow [T(Gf_t)]\sset T^1_\cG f_o$
  if the surjectivity holds for all the finite jets.
  Recall the general fact: the map $f_t:G\to T^1_\cG f_o$ is separable (as a morphism onto its image) iff $T_G f_t\twoheadrightarrow[T Gf_t]$.
   (Namely, the surjectivity holds for all the finite jets.)
\bed The unfolding
$F$ is called separable  if the map $G\stackrel{f_t}{\to}T^1_\cG f_o$ is separable.
\eed
If a field $\k$ is not algebraically closed then we call $F$ separable if $\bar\k\otimes F$ is separable.
\bex
\bee[\bf i.]
\item For $char(\k)=0$ the group orbit map is always separable. Hence any unfolding is separable.
\\\parbox{11cm}{\item The (in)separability is preserved under the $\cG_t$-equivalence.
 For example, for $\cG=\cR$ and any $\cb$ we have the commutative diagram of algebraic varieties}
\quad $\bM jet_\cb(G)\stackrel{f_t}{\to}jet_\cb(T^1_\cR f_o)\\\downarrow f_t\circ\phi\ \ \|\\jet_\cb T^1_\cR (f_o\circ\phi)=jet_\cb T^1_\cR (f_o)\eM$.
 \item A trivial unfolding is separable. Indeed, we can take $F=(f_o,t)$, then $G\stackrel{f_o}{\to}T^1_\cG f_o$ is the zero map.
  Hence $[T_G f_t]=0$. And thus trivially $T_G f_t\twoheadrightarrow [T(G f_t)]$.
\eee
\eex
\bel
An unfolding $F$ is inseparable iff $f\stackrel{\cG_t}{\sim}f_o+t^df_d+\cdots$, where $char(\k)\mid d$ and $f_d\not\in T_\cG f_o$.
\eel
\bpr
$\Lleftarrow$ We can assume $f=f_o+t^df_d+\cdots$ with  $char(\k)\mid d$ and $f_d\not\in T_\cG f_o$.
 Then $T_G f_t\sseteq (t)^{d+1}\cdot \RpXt$.
 Thus for $\cb=(t^{d+1},(x)^N)$ we get $jet_\cb(T_G f_t)=0$.
 But $G f_t\ni f_o+ c^d t^d f_d+\cdots$ for any $c\in \k$. Therefore the subvariety $[jet_\cb(G f_t)]\sset jet_\cb T^1_\cG f_o$ is of positive dimension. Thus
   $ T(jet_\cb G f_t)\neq0$. Therefore the map $T_G f_t\to [TG(f_t)]\sset T^1_\cG f_o$ cannot be surjective.

$\Rrightarrow$ Suppose  $f\sim f_o+t^df_d+\cdots$ with  $char(\k)\nmid d$ and $f_d\not\in T_\cG f_o$.
Then $jet_\cb T_G f_t=Span_{\k[[t]]}\{t^{d-1}f_d+\cdots\}$
 and $T jet_\cb(G f_t)=Span_{\k[[t]]}\{t^{d-1}f_d+\cdots\}$. Hence $F$ is separable.

Otherwise for each $d$ we get: $f\sim f_o+t^df_d+\cdots$ with $f_d\in T_\cG f_o$. Then $F$ is formally trivial. Trivialize it to get: $(f_o,t)$ is separable.
\epr

\end{document}